# Topological Quantum Computation from the 3-dimensional Bordism 2-Category

Juan Ospina
Logic and Computation Group
Mathematics Engineering Program
Department of Mathematics
School of Sciences
EAFIT University
Medellín, Colombia
jospina@eafit.edu.co

**ABSTRACT**

A great part of the mathematical foundations of topological quantum computation is given by the theory of modular categories which provides a description of the topological phases of matter such as anyon systems. In the near future the anyonic engineering will provide the anyonic devices from which the topological quantum computers will be constructed. From other side the string anyons are interesting topological phases of matter which can be described using mathematical constructions such as Frobenius algebras and open-closed string topological quantum field theories which are based on cobordism categories. Recently was proposed that is possible to obtain representations of cobordism categories using modular categories. In the present work, the modular categories resulting as representations of the 3-dimensional bordism 2-category are used with the aim to construct a new model of topological quantum computation. Such new model is named "Sanyon Topological Quantum Computation" and it is is theoretically performed by evolving non-abelian string- anyons (sanyons) using the Loop Braid Group and the open-closed cobordism category. The output of the computation uniquely depends on how the sanyons have been braided by the Loop Braid Group and operated by the generators of the cobordism category. Small disturbances do not unravel the loop braids and the cobordisms, making the computation resistant to errors and decoherence.

**Keywords:**. Topological Quantum Computation, Frobenius Algebras, Open-Closed Cobordism Category, Loop Braid Group, String anyons, Automated Theorem Proving, Quantum Automated Theorem Provers.

## 1. INTRODUCTION

In a topological quantum computer [1,2] the quantum information is encoded and processed using the so called topological quantum bits. The Hilbert space associated with the quantum states of a topological qubit is a non-local functional structure and for this reason the topological quantum bits are intrinsically protected against the decoherence. One possible physical implementation of a topological qubit is reached using the so called Majorana fermion zero modes corresponding to superconducting nanowires with strong spin-orbit couplings. Recently was proposed that a network of coupled superconducting nanowires with Majorana fermion zero modes is able to support a powerful class of topological non-Abelian particle named genon [3,4]. The braiding of such genons provide topological quantum gates which are universal.

The standard genom is a particle without extension. In this paper we consider a string genon or sgenon for which the braiding is mathematically performed using the Loop Braid Group [5,6,7] and the open-closed cobordism category [8,9].

The open-closed cobordism category is studied using the notion of Frobenius algebra. In the next section the basic facts about monoidal categories, Frobenius algebras and their implementation using automatic theorem provers are presented. As is pointed in [10], it is expected that automated theorem provers (ATPs) will play an important role in quantum computing in general and in topological quantum computing in particular; given the crucial role that the automated theorem provers are playing actually in classical computing. Specifically we will use automatic theorem provers which admit the THO language based on typed higher order logic. Some examples of such ATPs are: `agsyHOL, Isabelle-HOT, Leo II, Satallax and cocATP`. In this work `Leo II`, will be used extensively [11].

At the section fourth we will use all the code generated at section 2 with the aim to automatically prove theorems in the open-closed cobordism category. In the third section a quantum model for Khovanov homology for tangles will be presented. The string-anyon topogological quantum computer which will be structured using the automatic proved mathematics at sections 2 and 4 , is able to compute the topological invariant to be presented at the third section.

## 2. AUTOMATIC FORMALIZATION OF MATHEMATICS FOR TOPOLOGICAL QUANTUM COMPUTATION

A weak monoidal category satisfies the following pentagon identity for all objects *x* , *y*, *z*, *w* in the category [8,12]

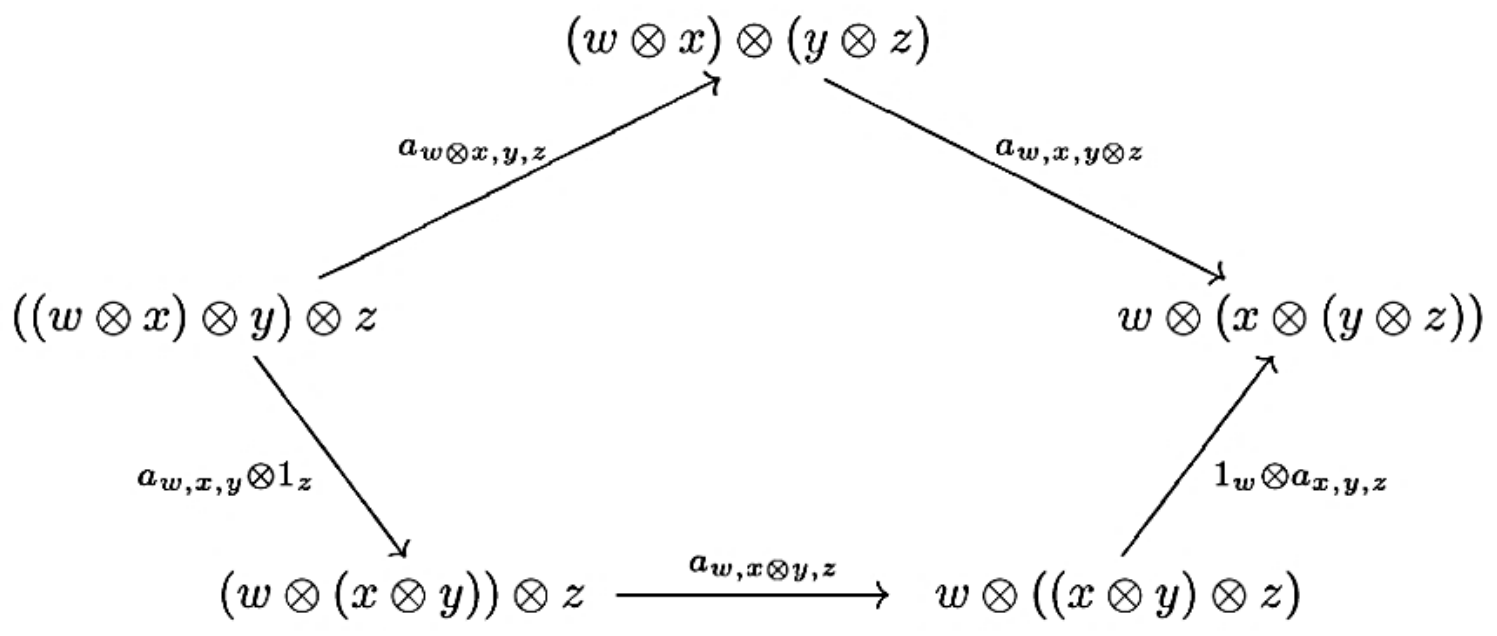

(0.0)

where

$$a_{x,y,z}: (x \otimes y) \otimes z \to x \otimes (y \otimes z) \tag{0.0.A}$$

is one natural isomorphism.

The pentagon identity (0.0) is proved using the automatic theorem prover `Leo II` with the following commands written with TPTP thf language [11]:

```
thf(c_type,type,(
    c: $tType )).
thf(one,type,(
    one: c )).
thf(alpha_decl,type,(alpha: c > c  )).
thf(iden_decl,type,(iden: c > c  )).
thf(multo_decl,type,(multo: c > c > c )).
thf(multm_decl,type,(multm: (c > c) > (c > c ) > (c > c ) )).
thf(axio1,axiom,(! [X: c, Y: c, Z: c] :
       ( (alpha @ (multo @ (multo @ X @ Y) @ Z))  = (multo @ X @ (multo @ Y @ Z))
)    )).
thf(axio2,axiom,(! [X: c] :
       ( (iden @ X)  = X     )    )).
thf(axio3,axiom,(! [X: c, Y: c, Z:c, W: c] :
```

```
      ( ((multm @ alpha @ iden) @ (multo @ (multo @ (multo @ W @ X ) @ Y) @ Z))
= (multo @ (multo @ W @ (multo @ X @ Y ) ) @ Z)      )     )).
thf(axio3A,axiom,(! [X: c, Y: c,Z:c,W:c] :
      ( ((multm @ iden @ alpha) @ (multo @ W @ (multo @ (multo @ X @ Y ) @ Z) )
)   = (multo @ W @ (multo @ X @ (multo @ Y @ Z)))       )     )).
thf(conje,conjecture,(! [X:c, Y:c, Z:c, W:c] :
      (  ((multm @ iden @ alpha) @ (alpha @ ((multm @ alpha @ iden) @ (multo @
(multo @ (multo @ W @ X) @ Y) @ Z) )) )   =  (alpha @ (alpha @ (multo @ (multo @
(multo @ W @ X ) @ Y) @ Z)))          )     )).
```

The corresponding output generated by `LeoII` is

```
% END OF SYSTEM OUTPUT
% RESULT: SOT_l9PD5B - LEO-II---1.6.2 says Theorem - CPU = 0.01 WC = 0.04
% OUTPUT: SOT_l9PD5B - LEO-II---1.6.2 says CNFRefutation - CPU = 0.01 WC = 0.04
```

A weak monoidal category satisfies the following triangle identity for all objects $x$ , $y$ in the category

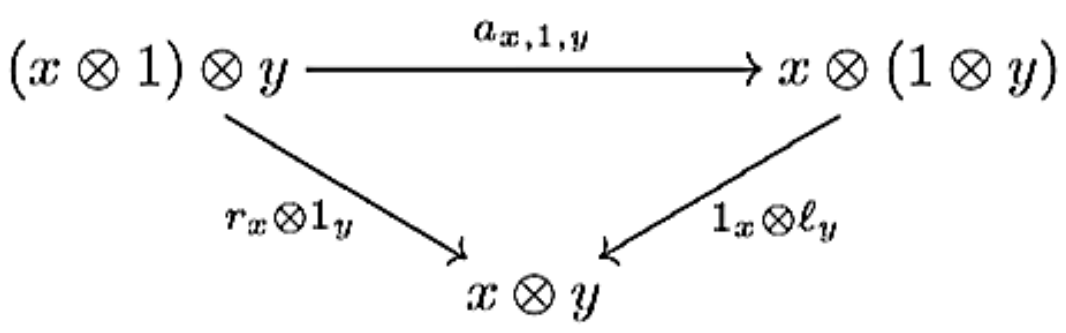

(0.0.B)

where

$$\ell_x: 1 \otimes x \to x,$$

$$r_x: x \otimes 1 \to x, \qquad (0.0C)$$

are other two natural isomorphisms.

The identity (0.0B) is proved adding the following commands to the code previously presented

```
thf(left_decl,type,(left: c > c  )).
thf(right_decl,type,(right: c > c  )).
thf(axio2A,axiom,(! [X: c] :
      ( (left @ (multo @ one @ X))  = X      )     )).
thf(axio2B,axiom,(! [X: c] :
      ( (right @ (multo @ X @ one))  = X      )     )).
thf(axio2C,axiom,(! [X: c, Y: c] :
      ( ((multm @ right @ iden) @ (multo @ (multo @ X @ one) @ Y))  = (multo @ X
@ Y)      )     )).
thf(axio2D,axiom,(! [X: c, Y:c] :
      ( ((multm @ iden @ left) @ (multo @ X @ (multo @ one @ Y) ))  = (multo @ X
@ Y)      )     )).
thf(conje2,conjecture,(! [X:c, Y:c] :
      ( ((multm @ iden @ left) @ (alpha @ (multo @ (multo @ X @ one) @ Y)) ) =
((multm @ right @ iden) @ (multo @ (multo @ X @ one ) @  Y))
)     )).
```

The corresponding output generated by `LeoII` is

```
% END OF SYSTEM OUTPUT
```

```
% RESULT: SOT_KytOVn - LEO-II---1.6.2 says Theorem - CPU = 0.02 WC = 0.06
% OUTPUT: SOT_KytOVn - LEO-II---1.6.2 says CNFRefutation - CPU = 0.02 WC = 0.07
```

The Frobenius algebras in a generic symmetric monoidal category satisfy the so called pentagon identity given by

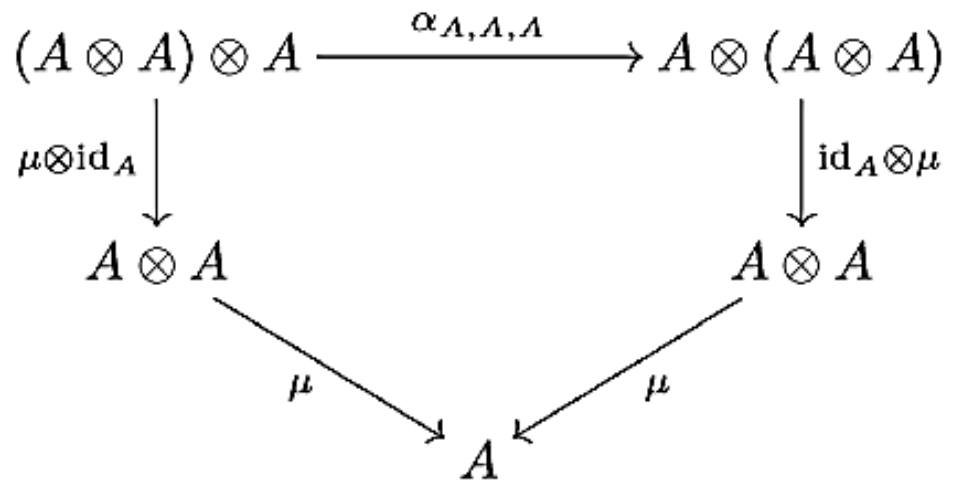

(0.1)

or equivalently by

$\mu\,(\mu \otimes I)\,|(A \otimes A) \otimes A\rangle = \mu\,(I \otimes \mu)\,\alpha|(A \otimes A) \otimes A\rangle$ (0.2)

The pentagon identity (0.1) or (0.2) is codified in the automatic theorem prover Leo II using the following commands

```
thf(axio3,axiom,(! [X: $aaxa] :
      ( (mu @ (idmu @ (alpha @ X))   )  = (mu @ (muid @ X) )   )    )).
```

We are using the TPTP thf language with the following specifications

```
thf(alpha_decl,type,(alpha: $aaxa > $axaa)).
thf(invalpha_decl,type,(invalpha: $axaa > $aaxa )).
thf(mu_decl,type,(mu: $aa > $a )).
thf(muid_decl,type,(muid: $aaxa > $aa )).
thf(idmu_decl,type,(idmu: $axaa > $aa )).
thf(id_decl,type,(id: $a > $a )).
```

The Frobenius algebras in a generic symmetric monoidal category satisfy the so called co-pentagon identity given by

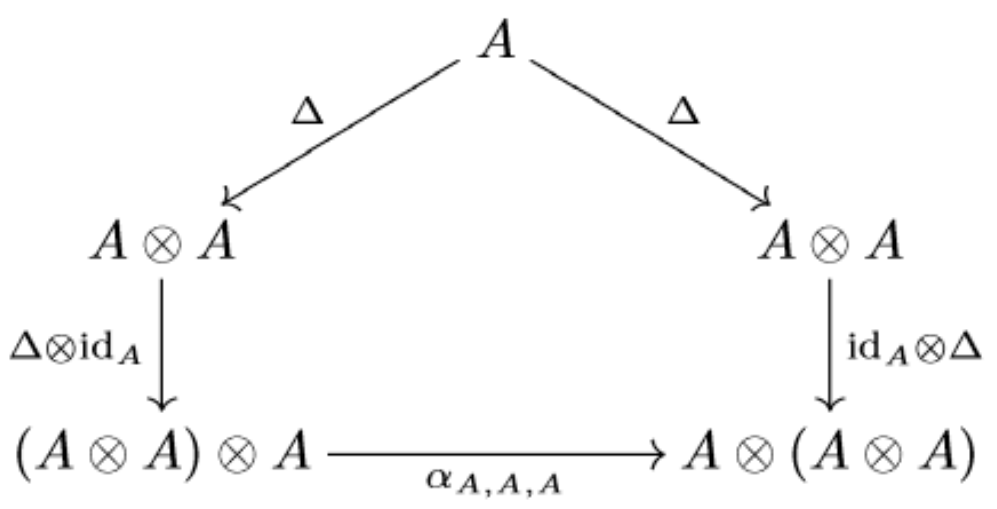

(0.3)

or equivalently by

$\alpha\,(\Delta \otimes I)\Delta\,|A\rangle = (I \otimes \Delta)\,\Delta|A\rangle$ (0.4)

The pentagon identity (0.3) or (0.4) is codified in the automatic theorem prover Leo II using the following commands

```
thf(axio4,axiom,(! [X: $a] :
      ( (iddelta @ (delta @ X) )  = (alpha @ (deltaid @ (delta @ X)) ) )   )).
```

We are using the TPTP thf language with the following specifications

```
thf(delta_decl,type,(delta: $a > $aa )).
thf(deltaid_decl,type,(deltaid: $aa > $aaxa )).
thf(iddelta_decl,type,(iddelta: $aa > $axaa )).
```

The Frobenius algebras in a generic symmetric monoidal category satisfy the so called triangle identity given by

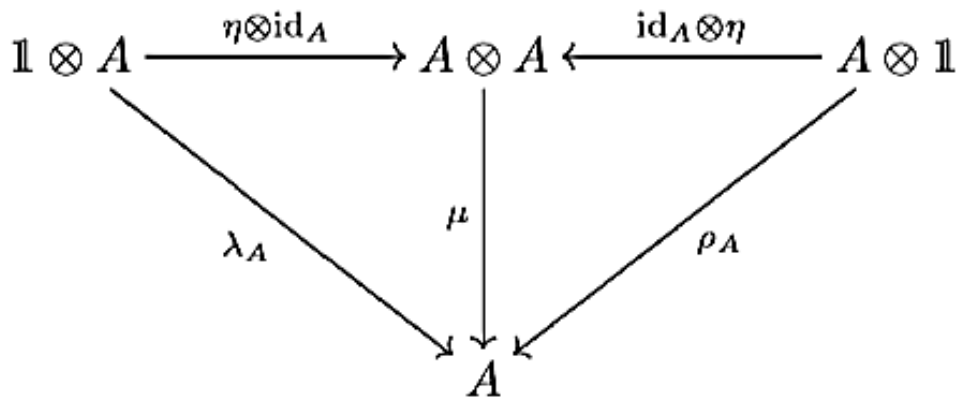

(0.5)

or equivalently by

$\mu\,(\eta \otimes I)\,(|1\rangle \otimes |A\rangle\,) \;=\; \lambda\,(|1\rangle \otimes |A\rangle\,)$ (0.6)

$\mu\,(I \otimes \eta)\,(|A\rangle \otimes |1\rangle\,) \;=\; \rho\,(|A\rangle \otimes |1\rangle\,)$ (0.7)

The triangle identity (0.5) or (0.6)-(07) is codified in the automatic theorem prover Leo II using the following commands

```
thf(axio1,axiom,(! [X: $ja] :
      ( (lamb @ X)  = (mu @ (etaid @ X ) )   )    )).
thf(axio2,axiom,(! [X: $ai] :
      ( (rho @ X)  = (mu @ (ideta @ X ) )   )    )).
```

We are using the TPTP thf language with the following specifications

```
thf(eta_decl,type,(eta: $i > $a )).
thf(etaid_decl,type,(etaid: $ja > $aa )).
thf(ideta_decl,type,(ideta: $ai > $aa )).
thf(lamb_decl,type,(lamb: $ja > $a )).
thf(rho_decl,type,(rho: $ai > $a )).
```

The Frobenius algebras in a generic symmetric monoidal category satisfy the so called co-triangle identity given by

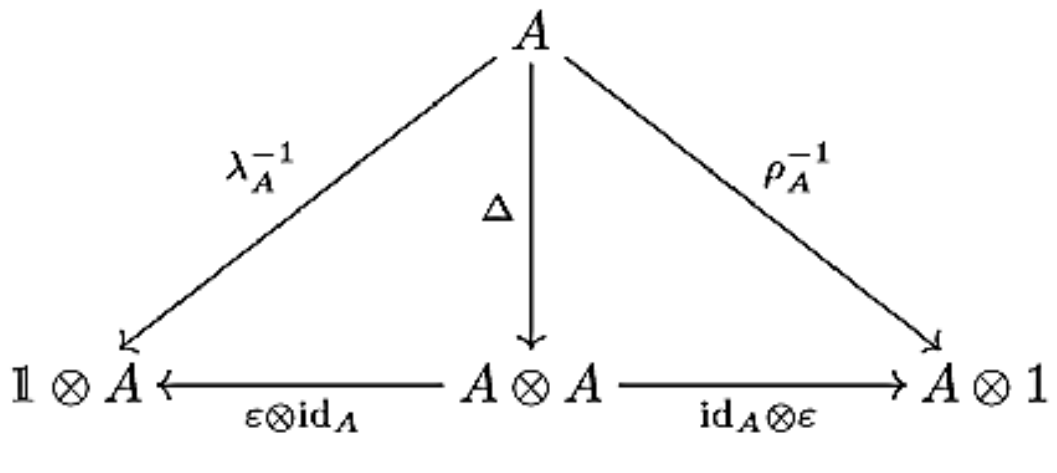

(0.8)

or equivalently by

$(\varepsilon \otimes I)\Delta\,(|A\rangle\,) \;=\; \lambda^{-1}\,(\,|A\rangle\,)$ (0.9)

$(I\otimes\varepsilon)\,\Delta\,(|A\rangle) \;=\; \rho^{-1}\,(|A\rangle)$ (0.10)

The co-triangle identity (0..8) or (0.9)-(0.10) is codified in the automatic theorem prover Leo II using the following commands

```
thf(epsilon_decl,type,(epsilon: $a > $i )).
thf(invlamb_decl,type,(invlamb: $a > $ja )).
thf(invrho_decl,type,(invrho: $a > $ai )).
thf(epsilonid_decl,type,(epsilonid: $aa > $ja )).
thf(idepsilon_decl,type,(idepsilon: $aa > $ai )).
thf(axio5,axiom,(! [X: $a] :
       (  (epsilonid @ (delta @ X) )  = (invlamb @ X )     )    )).
thf(axio6,axiom,(! [X: $a] :
       (  (idepsilon @ (delta @ X) )  = (invrho @ X )     )    )).
```

Now, in a braided monoidal category the following equation is satisfied:

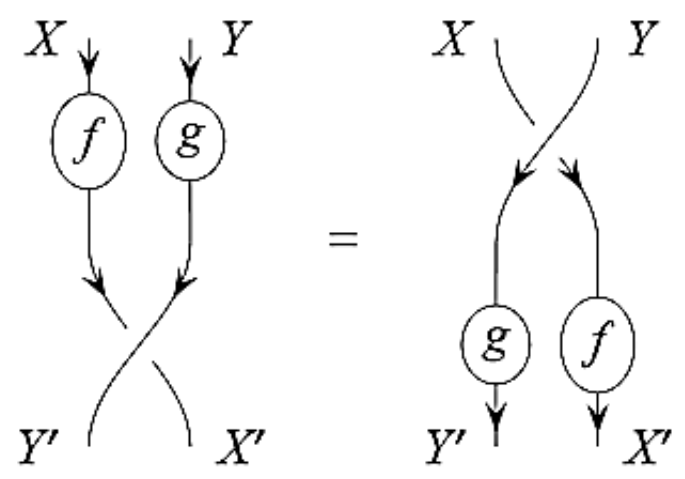

(0.11)

or equivalently

$$\sigma\, ((f \otimes g)\, (|X\rangle \otimes |Y\rangle\, ))\;\; = (g\otimes f)\;\; (\sigma\, (|X\rangle \otimes |Y\rangle\, )\, ) \qquad (0.12)$$

The equation (0.11) or (0.12) is proved simultaneously by Isabelle-HOT, Leo II and Satallax using the following code

```
thf(c_type,type,(c: $tType )).
thf(cc_type,type,(cc: $tType )).
thf(braid_decl,type,(braid: cc > cc  )).
thf(up_decl,type,(up: c > c  )).
thf(down_decl,type,(down: c > c  )).
thf(multo_decl,type,(multo: c > c > cc  )).
thf(multm_decl,type,(multm: (c > c) > (c > c ) > (cc > cc ) )).
thf(axio,axiom,(![X:c, Y:c, F: c > c, G: c > c]:
       ( ((multm @ F @ G) @ (multo @ X @ Y))  = (multo @ (F @ X) @ (G @ Y))
))).
thf(axio1,axiom,(![X:c, Y:c]:
       (  (braid  @  (multo  @  X  @  Y))    =  (multo  @  (up  @  Y)  @  (down  @  X))
))).
thf(axio1A,axiom,(![X:c, F:c > c]:
       ( (F @ (up @ X))  = (up @ (F @ X))          ))).
thf(axio1B,axiom,(![X:c, F:c > c]:
       ( (F @ (down @ X))  = (down @ (F @ X))          ))).
thf(conje,conjecture,(! [X:c, Y:c, F:c > c, G: c > c] :
       ( (braid @ ( (multm @ F @ G) @ (multo @ X @ Y)  )) =   ( (multm @ G @ F) @
(braid @ (multo @ X @ Y))    )
)    )).
```

and the corresponding output from the ATPs is

```
% END OF SYSTEM OUTPUT
% RESULT: SOT_Bbh0wG - Isabelle-HOT---2015 says Theorem - CPU = 6.25 WC = 5.36
SolvedBy = simp
% OUTPUT: SOT_Bbh0wG - Isabelle-HOT---2015 says Assurance - CPU = 6.25 WC = 5.36

% END OF SYSTEM OUTPUT
% RESULT: SOT_Bbh0wG - LEO-II---1.6.2 says Theorem - CPU = 0.00 WC = 0.02
% OUTPUT: SOT_Bbh0wG - LEO-II---1.6.2 says CNFRefutation - CPU = 0.00 WC = 0.02

% END OF SYSTEM OUTPUT
% RESULT: SOT_Bbh0wG - Satallax---2.8 says Theorem - CPU = 0.00 WC = 0.15
% OUTPUT: SOT_Bbh0wG - Satallax---2.8 says Proof - CPU = 0.00 WC = 0.16
```

The corresponding directed acyclic graph generated using `IDV` from the proof given by `Leo II` is showed at the following figure

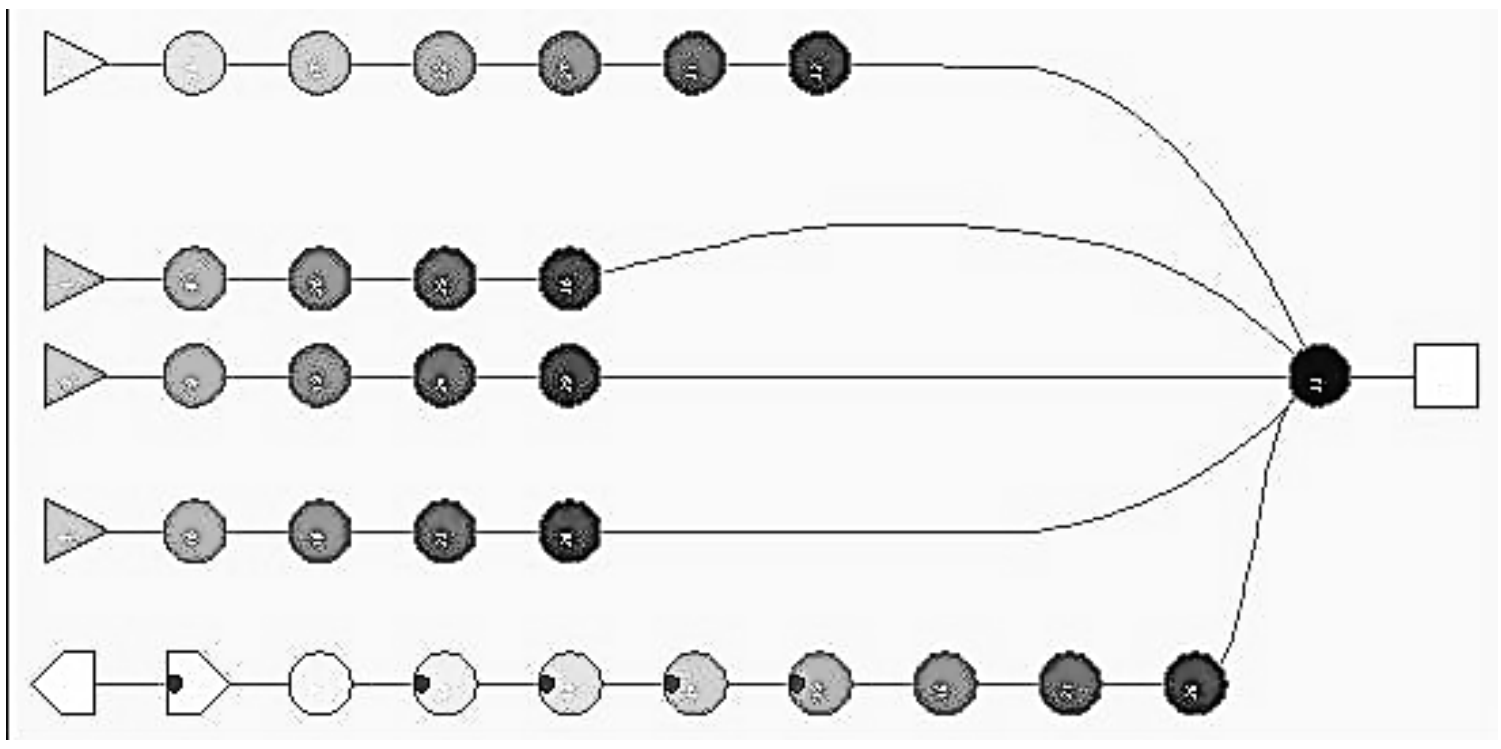

Also, in a braided monoidal category the following equation is satisfied

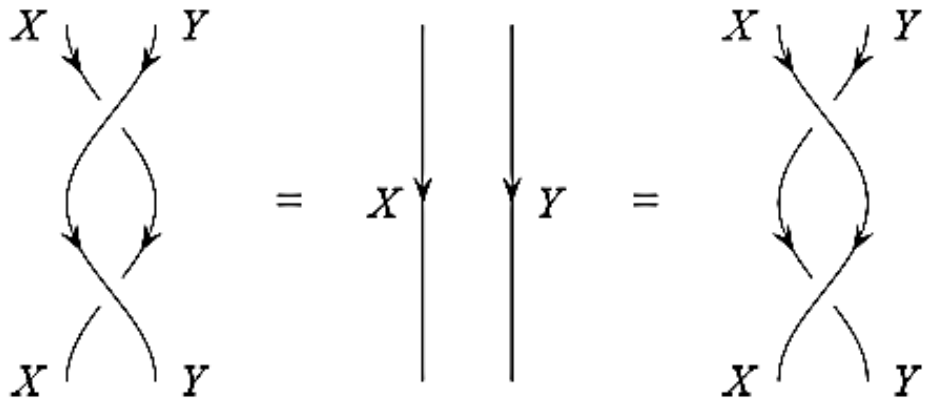

(0.13)

or equivalently

$$\sigma^{-1}(\sigma|(X\rangle \otimes |Y\rangle )) = \sigma (\sigma^{-1}(|X\rangle \otimes |Y\rangle )) = (I |X\rangle) \otimes (I |Y\rangle) = (|X\rangle \otimes |Y\rangle. \quad (0.14)$$

The equation (0.13) or (0.14) is proved simultaneously by Isabelle-HOT, Leo II using the following code

```
thf(c_type,type,(c: $tType )).
thf(cc_type,type,(cc: $tType )).
thf(braid_decl,type,(braid: cc > cc  )).
thf(up_decl,type,(up: c > c  )).
thf(down_decl,type,(down: c > c  )).
thf(invbraid_decl,type,(invbraid: cc > cc  )).
thf(iden_decl,type,(iden: c > c  )).
```

```
thf(multo_decl,type,(multo: c > c > cc  )).
thf(multm_decl,type,(multm: (c > c) > (c > c ) > (cc > cc ) )).
thf(axio0,axiom,(![X:c]:( (iden @ X)  = (X)          ))).
thf(axio,axiom,(![X:c, Y:c, F: c > c, G: c > c]:
       ( ((multm @ F @ G) @ (multo @ X @ Y))  = (multo @ (F @ X) @ (G @ Y))
))).
thf(axio1,axiom,(![X:c, Y:c]:
       ( (braid @ (multo @ X @ Y))  = (multo @ (up @ Y) @ (down @ X))
))).
thf(axio1A,axiom,(![X:c, Y:c]:
       ( (invbraid @ (multo @ X @ Y))  = (multo @ (down @ Y) @ (up @ X))
))).
thf(axio2,axiom,(![X:c]:
       ( (up @ (up @ X ))  = (X)          ))).
thf(axio2A,axiom,(![X:c]:
       ( (down @ (down @ X ))  = (X)          ))).
thf(conje,conjecture,(! [X:c, Y:c] :
       (   (invbraid @ (braid @ (multo @ X @ Y)) ) = (braid @ (invbraid @ (multo
@ X @ Y)))
)    )).
thf(conje1,conjecture,(! [X:c, Y:c] :
       (   (invbraid @ (braid @ (multo @ X @ Y)) ) = ((multm @ iden @ iden) @
(multo @ X @ Y))
)    )).
thf(conje2,conjecture,(! [X:c, Y:c] :
       (   (braid @ (invbraid @ (multo @ X @ Y))) = ((multm @ iden @ iden) @
(multo @ X @ Y))
)    )).
```

In a braided monoidal category the following hexagon identities are valid:

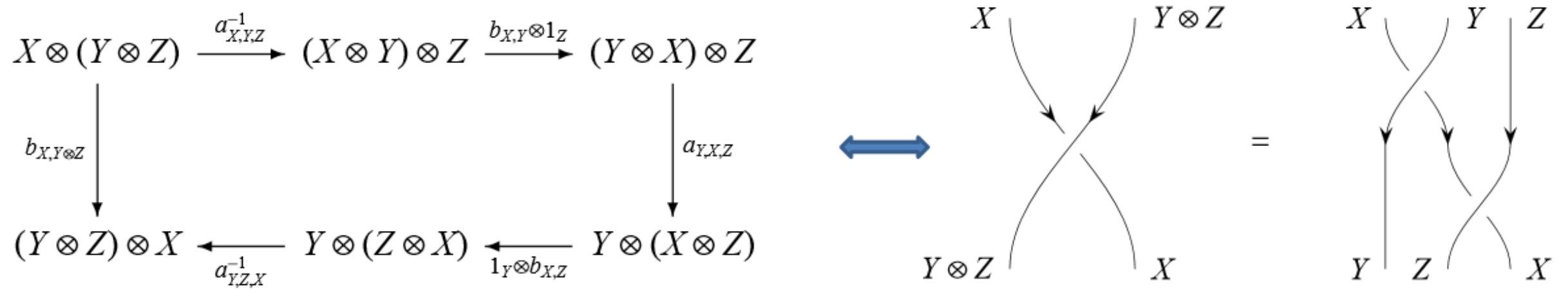

(0.15)

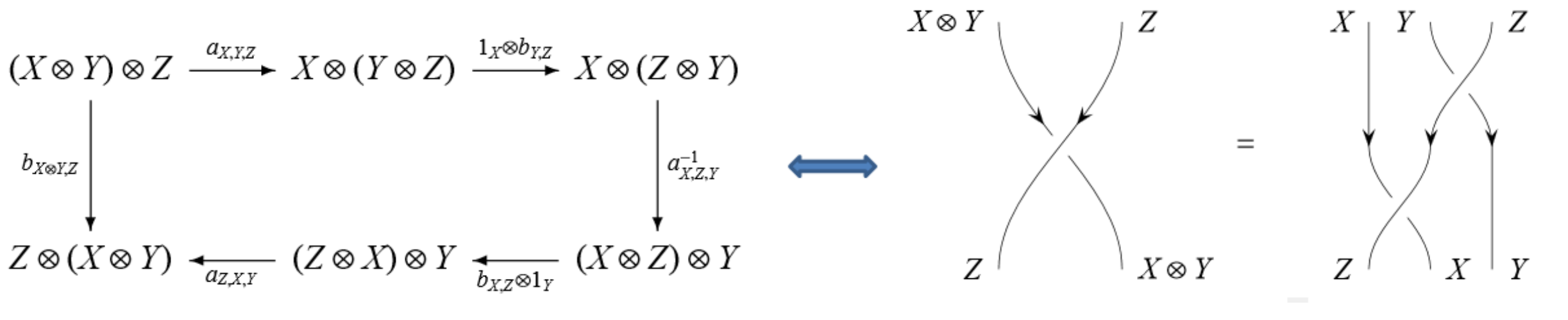

(0.16)

The equation (0.15) and (0.16) are formalized simultaneously by Isabelle-HOT, Leo II using the following code

```
thf(c_type,type,(c: $tType )).
thf(cc_type,type,(cc: $tType )).
thf(ccxc_type,type,(ccxc: $tType )).
thf(cxcc_type,type,(cxcc: $tType )).
```

```
thf(alpha_decl,type,(alpha: ccxc > cxcc  )).
thf(invalpha_decl,type,(invalpha: cxcc > ccxc  )).
thf(braid_decl,type,(braid: cc > cc  )).
thf(up_decl,type,(up: c > c  )).
thf(down_decl,type,(down: c > c  )).
thf(up1_decl,type,(up1: cc > cc  )).
thf(down1_decl,type,(down1: cc > cc  )).
thf(invbraid_decl,type,(invbraid: cc > cc  )).
thf(braid1_decl,type,(braid1: cxcc > ccxc  )).
thf(braid2_decl,type,(braid2: ccxc > cxcc  )).
thf(iden_decl,type,(iden: c > c  )).
thf(multo_decl,type,(multo: c > c > cc  )).
thf(multo1_decl,type,(multo1: cc > c > ccxc  )).
thf(multo2_decl,type,(multo2: c > cc > cxcc  )).
thf(multm_decl,type,(multm: (c > c) > (c > c ) > (cc > cc ) )).
thf(multm1_decl,type,(multm1: (cc > cc) > (c > c ) > (ccxc > ccxc ) )).
thf(multm2_decl,type,(multm2: (c > c) > (cc > cc ) > (cxcc > cxcc ) )).
thf(axio0,axiom,(![X:c]:( (iden @ X)  = (X)          ))).
thf(axio0A,axiom,(![X:c]:( (down @ (down @ X ))  = (X)          ))).
thf(axio0B,axiom,(![X:c]:( (up @ (up @ X ))  = (X)          ))).
thf(axio,axiom,(![X:c, Y:c, F: c > c, G: c > c]:
      ( ((multm @ F @ G) @ (multo @ X @ Y))  = (multo @ (F @ X) @ (G @ Y))
))).
thf(axioA,axiom,(![X:c, Y:c,Z:c, F: cc > cc, G: c > c]:
      ( ((multm1 @ F @ G) @ (multo1 @ (multo @ X @ Y) @ Z))  = (multo1 @ (F @
(multo @ X @ Y)) @ (G @ Z))          ))).
thf(axioAA,axiom,(![X:c, Y:c,Z:c, F: c > c, G: cc > cc]:
      ( ((multm2 @ F @ G) @ (multo2 @ X @ (multo @ Y @ Z) )  )  = (multo2 @ (F @
X) @ (G @ (multo @ Y @ Z))   )          ))).
thf(axio1,axiom,(![X:c, Y:c]:
      ( (braid @ (multo @ X @ Y))   = (multo @ (up @ Y) @ (down @ X))
))).
thf(axio1A,axiom,(![X:c, Y:c]:
      ( (invbraid @ (multo @ X @ Y))   = (multo @ (down @ Y) @ (up @ X))
))).
thf(axio1B,axiom,(![X:c, Y:c,Z:c]:
      ( (braid1 @ (multo2 @ X @ (multo @ Y @ Z)))  = (multo1 @ (up1 @ (multo @ Y
@ Z )) @ (down @X))          ))).
thf(axio1C,axiom,(![X:c, Y:c,Z:c]:
      ( (braid2 @ (multo1 @ (multo @ X @ Y) @ Z  ))  = (multo2 @ (up @ Z) @
(down1 @ (multo @ X @ Y )) )          ))).
thf(axio2,axiom,(![X:c, Y:c, Z:c]:
      ( (alpha @ (multo1 @ (multo @ X @ Y) @ Z))  = (multo2 @ X @ (multo @ Y @ Z
))          ))).
thf(axio2A,axiom,(![X:c, Y:c, Z:c]:
     ( (invalpha @ (multo2 @ X @ (multo @ Y @ Z)) )  = (multo1 @ (multo @ X @ Y)
@ Z)          ))).
thf(conje,conjecture,(! [X:c, Y:c, Z:c] :
      (    ( invalpha @ ((multm2 @ iden @ braid) @ (alpha @((multm1 @ braid @
iden)  @ (invalpha @ (multo2 @ X @ (multo @ Y @ Z))) ))  )) = (multo1 @ (multo @
(up @ Y) @ (up @ Z ) ) @ X)
)    )).
thf(conje1,conjecture,(! [X:c, Y:c, Z:c] :
```

```
       (     ( alpha @  ((multm1  @  braid  @  iden)  @  (invalpha  @((multm2  @  iden  @
braid)   @  (alpha  @  (multo1  @  (multo  @  X  @  Y)  @  Z  ))  ))    ))  =  (multo2  @  Z  @
(multo @  (down  @  X)  @  (down  @  Y)  ))
)     )).
```

With this formalization the hexagon identities (0.15) and (0.16) are reduced to the following axioms

```
thf(axioH1,axiom,(![X:c, Y:c, Z:c]:
       (  (invalpha  @  ((multm2  @  iden  @  braid)  @  (alpha  @((multm1  @  braid  @  iden)
@  (invalpha  @  (multo2  @  X  @  (multo  @  Y  @  Z)))  ))   ))   =  (braid1  @  (multo2  @  X  @
(multo  @  Y  @  Z))  )          ))).
thf(axioH2,axiom,(![X:c, Y:c, Z:c]:
       (  (alpha  @  ((multm1  @  braid  @  iden)  @  (invalpha  @((multm2  @  iden  @  braid)
@  (alpha  @  (multo1  @  (multo  @  X  @  Y)  @  Z  ))  ))   ))   =  (braid2  @  (multo1  @  (multo
@  X  @  Y)  @  Z   ))          ))).
```

Finally, in a braided monoidal category it is possible to prove the Yang-Baxter equation

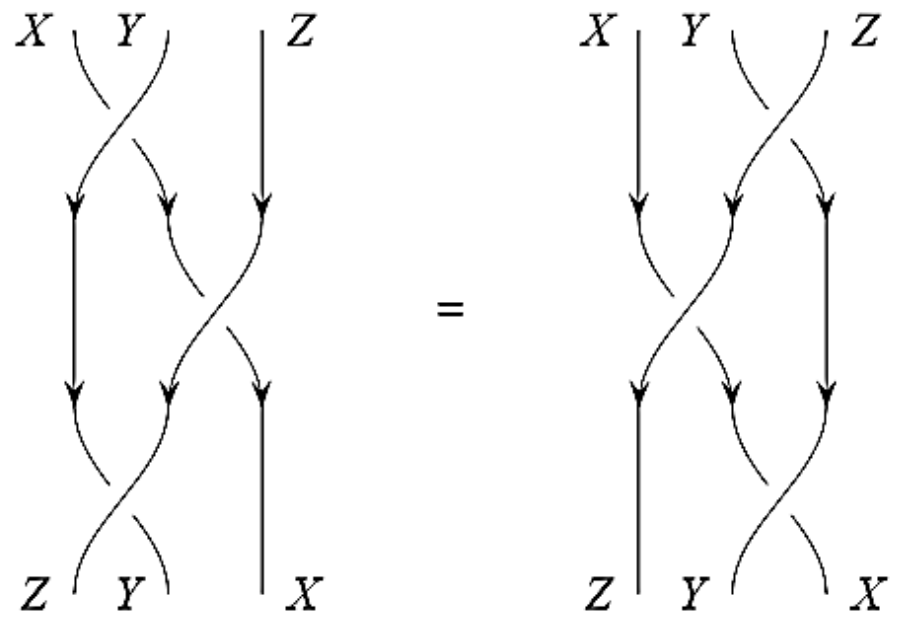

(0.17)

using the following code

```
thf(c_type,type,(c: $tType )).
thf(cc_type,type,(cc: $tType )).
thf(ccxc_type,type,(ccxc: $tType )).
thf(cxcc_type,type,(cxcc: $tType )).
thf(alpha_decl,type,(alpha: ccxc > cxcc  )).
thf(invalpha_decl,type,(invalpha: cxcc > ccxc  )).
thf(braid_decl,type,(braid: cc > cc  )).
thf(up_decl,type,(up: c > c  )).
thf(down_decl,type,(down: c > c  )).
thf(up1_decl,type,(up1: cc > cc  )).
thf(down1_decl,type,(down1: cc > cc  )).
thf(invbraid_decl,type,(invbraid: cc > cc  )).
thf(braid1_decl,type,(braid1: cxcc > ccxc  )).
thf(braid2_decl,type,(braid2: ccxc > cxcc  )).
thf(iden_decl,type,(iden: c > c  )).
thf(multo_decl,type,(multo: c > c > cc  )).
thf(multo1_decl,type,(multo1: cc > c > ccxc  )).
thf(multo2_decl,type,(multo2: c > cc > cxcc  )).
thf(multm_decl,type,(multm: (c > c) > (c > c ) > (cc > cc ) )).
thf(multm1_decl,type,(multm1: (cc > cc) > (c > c ) > (ccxc > ccxc ) )).
thf(multm2_decl,type,(multm2: (c > c) > (cc > cc ) > (cxcc > cxcc ) )).
thf(axio0,axiom,(![X:c]:(  (iden @ X)   = (X)          ))).
thf(axio0A,axiom,(![X:c]:(  (down @  (down @ X ))  = (X)          ))).
thf(axio0B,axiom,(![X:c]:(  (up @  (up @ X ))  = (X)          ))).
```

```
thf(axio0C,axiom,(![X:c, F:c > c]:
       ( (F @ (up @ X))  = (up @ (F @ X))          ))).
thf(axio0D,axiom,(![X:c, F:c > c]:
       ( (F @ (down @ X))  = (down @ (F @ X))          ))).


thf(axio,axiom,(![X:c, Y:c, F: c > c, G: c > c]:
       ( ((multm @ F @ G) @ (multo @ X @ Y))  = (multo @ (F @ X) @ (G @ Y))
))).
thf(axioA,axiom,(![X:c, Y:c,Z:c, F: cc > cc, G: c > c]:
       ( ((multm1 @ F @ G) @ (multo1 @ (multo @ X @ Y) @ Z))  = (multo1 @ (F @
(multo @ X @ Y)) @ (G @ Z))          ))).
thf(axioAA,axiom,(![X:c, Y:c,Z:c, F: c > c, G: cc > cc]:
       ( ((multm2 @ F @ G) @ (multo2 @ X @ (multo @ Y @ Z) )  )  = (multo2 @ (F @
X) @ (G @ (multo @ Y @ Z))   )          ))).
thf(axio1,axiom,(![X:c, Y:c]:
       ( (braid @ (multo @ X @ Y))   = (multo @ (up @ Y) @ (down @ X))
))).
thf(axio1A,axiom,(![X:c, Y:c]:
       ( (invbraid @ (multo @ X @ Y))   = (multo @ (down @ Y) @ (up @ X))
))).
thf(axio1B,axiom,(![X:c, Y:c,Z:c]:
       ( (braid1 @ (multo2 @ X @ (multo @ Y @ Z)))  = (multo1 @ (up1 @ (multo @ Y
@ Z )) @ (down @X))          ))).
thf(axio1C,axiom,(![X:c, Y:c,Z:c]:
       ( (braid2 @ (multo1 @ (multo @ X @ Y) @ Z  ))  = (multo2 @ (up @ Z) @
(down1 @ (multo @ X @ Y )) )          ))).
thf(axio2,axiom,(![X:c, Y:c, Z:c]:
       ( (alpha @ (multo1 @ (multo @ X @ Y) @ Z))  = (multo2 @ X @ (multo @ Y @ Z
))          ))).
thf(axio2A,axiom,(![X:c, Y:c, Z:c]:
      ( (invalpha @ (multo2 @ X @ (multo @ Y @ Z)) )  = (multo1 @ (multo @ X @ Y)
@ Z)          ))).
thf(conjeYB,conjecture,(! [X:c, Y:c, Z:c] :
       (    ( (multm1 @ braid @ iden ) @ (invalpha @ ((multm2 @ iden @ braid) @
(alpha @ ((multm1 @ braid @ iden) @ (multo1 @ (multo @ X @ Y) @ Z)))  )   )    )
= (invalpha @ ( (multm2 @ iden @ braid) @ (alpha @ ( (multm1 @ braid @ iden) @
(invalpha @ ( (multm2 @ iden @ braid ) @ (multo2 @ X @ (multo @ Y @ Z)))    )
))  )   )
)    )).
```

The corresponding directed acyclic graph generated using IDV from the proof given by Leo II is showed at the following figure

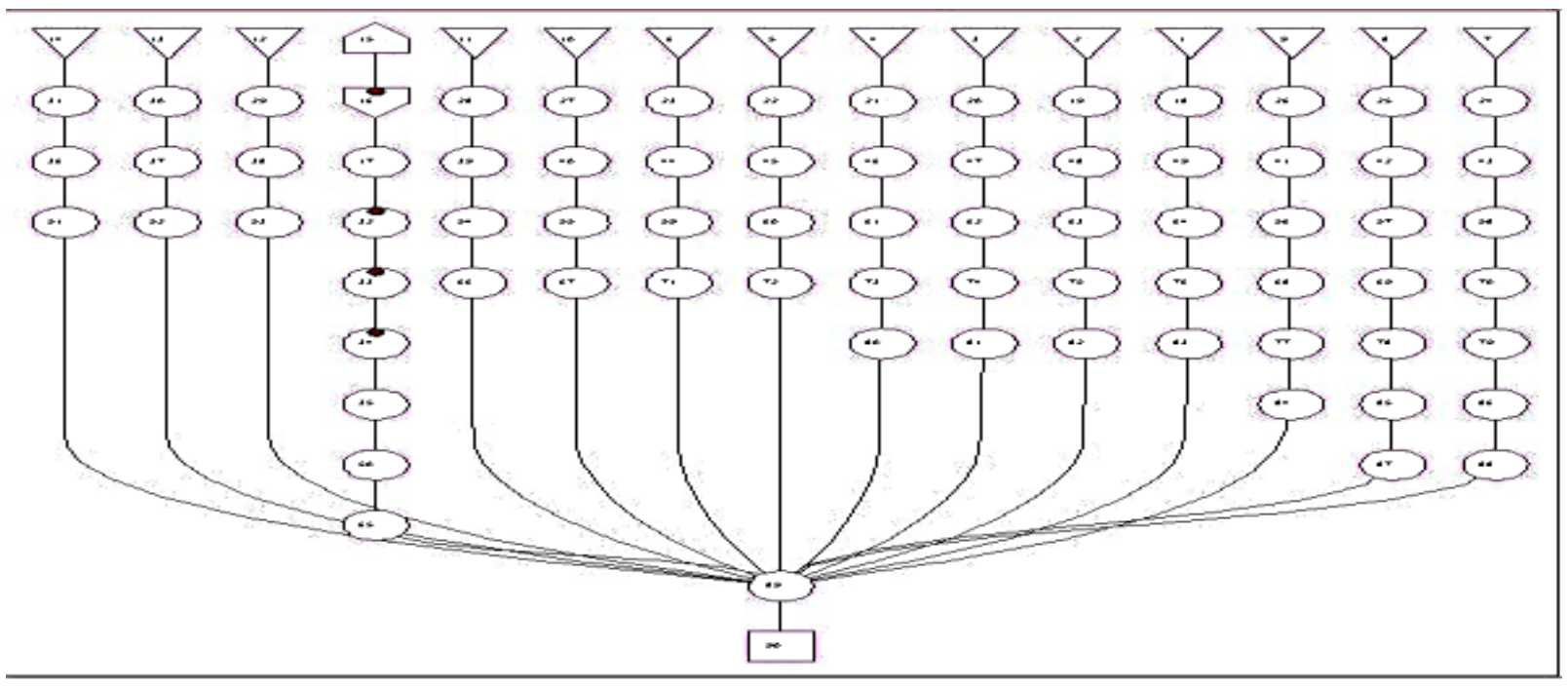

## 3. QUANTUM MODEL FOR KHOVANOV HOMOLOGY FOR TANGLES

A possible quantum model of Khovanov homology for tangles [2] will be build using super-symmetric string theory [13,14]. The $r$-th cohomology group of the tangle complex $C := [[T]]_{\varepsilon,Z}$ , is denoted by $H^r(C)$ . Assuming that $Z$ is an Euler-filtered open-closed String Topological Quantum Field Theory, then it is known that the filtration $F^*$ on the tangle complex $C$ induces a filtration on each homology group according to [8]

$$F^k\mathcal{H}^r(C) := \mathrm{image}\big(\mathcal{H}^r(F^kC) \hookrightarrow \mathcal{H}^r(C)\big) \tag{1}$$

where

$$F^kC := \bigoplus\nolimits_{r\in\mathbb{Z}} F^kC^r. \tag{2}$$

Such filtration defines a bigrading on the tangle homology given by

$$\mathcal{H}^{k,r}(C) := F^k\mathcal{H}^r(C)/F^{k+1}\mathcal{H}^r(C). \tag{3}$$

Assuming now, that $Z$ is an Euler-filtered open-closed String Topological Quantum Field Theory that satisfies Bar-Natan's conditions

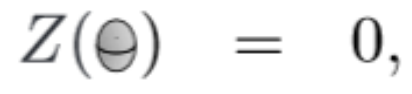

(4)

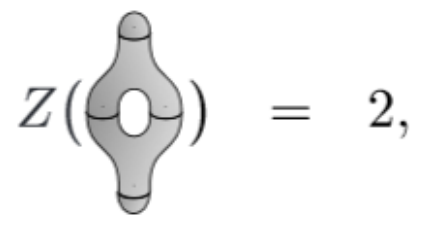

(5)

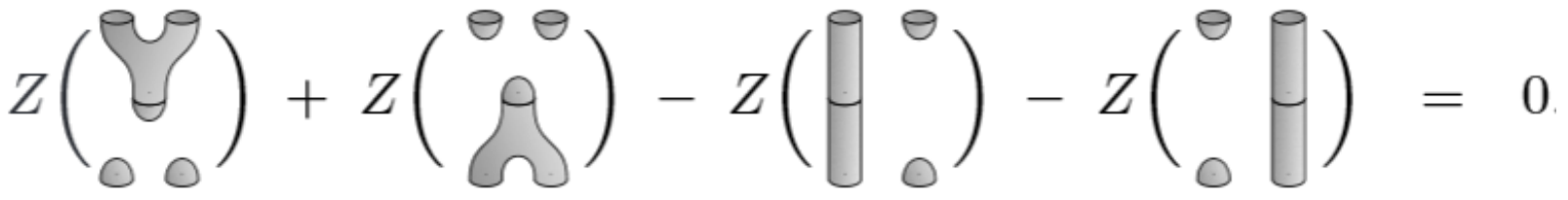

(6)

then it is known that the filtered Poincaré polynomial of the tangle complex $C$ given by

$$\mathcal{P}(T)_{\varepsilon,Z} := \sum_{r,k\in\mathbb{Z}} t^r A^k \mathrm{rk}\, \mathcal{H}^{k,r}(C), \tag{7}$$

is an invariant of the tangle $T$. Such invariant is called the 2-variable tangle polynomial and it satisfies

$$\mathcal{P}(T)_{\varepsilon,Z} \in \mathbb{Z}[t, t^{-1}, A, A^{-1}] \tag{8}$$

Now we find a quantum expression for (7) , which is able to mimic a string topological quantum computer for the computation of the tangle polynomial in Khovanov homology for tangles. We look for an expression which has the form of a generalized Witten index given by [13,14]

$$\lim_{\tau \to 0} \mathrm{Tr}(t^F A^p M \mathbf{e}^{(-\tau H)}) \tag{9}$$

The equation (9) is interpreted as the index of the Dirac-Ramond operator in the associated supersymmetric string theory. $P$ is a generalization of the world sheet momentum operator defined as

$$P = L_0 - \tilde{L}_0 \tag{10}$$

For the generalized Witten index given by (9) it is possible to write the following path integral

$$\lim_{\tau \to 0} \mathrm{Tr}(t^F A^p M \mathbf{e}^{(-\tau H)}) = \lim_{\tau \to 0} \int\int \mathbf{e}^{(-S(\phi, \psi))} \, dBoson \, dFermion \tag{11}$$

Now we study under which conditions (11) is a string topological quantum computer with the ability to compute the tangle polynomial (7) in the Khovanov homology for tangles. With this aim we introduce the following axioms which are connecting superstring theory with Khovanov homology for tangles.

**Axioms:** Given that $| rk >$ is an holological state for a given tangle in the Lauda formulation; or it is a supersymmetric state for a certain supersymmetric string quantum system in the Witten formulation, we introduce the following axioms for the quantum model of the Khovanov homology for tangles:

$$F \,| rk > = r \,| rk > \tag{12}$$

$$P \,| rk > = k \,| rk > \tag{13}$$

$$M \,| rk > = \mathrm{rank}(\, \mathcal{H}^{k,r}(C) \,) \,| rk > \tag{14}$$

$$U = t^F A^P M \tag{15}$$

where $U : C(T) \to C(T)$, is a unitary transformation in the Khovanov homology for the tangle $T$.

From the axioms that were introduced we derive the following propositions in Khovanov homology for tangles:

**Proposition 1.** $t^F \,| rk > = t^r \,| rk >$

**Proof:**

$$t^F \,| rk > = e^{\ln(t^F)} \,| rk >$$

$$= e^{(F \ln(t))} \,| rk >$$

$$= \left( \sum_{n=0}^{\infty} \frac{(F \ln(t))^n}{n!} \right) | rk >$$

$$= \sum_{n=0}^{\infty} \frac{\ln(t)^n F^n \mid rk>}{n!}$$

$$= \sum_{n=0}^{\infty} \frac{\ln(t)^n r^n \mid rk>}{n!}$$

$$= \left( \sum_{n=0}^{\infty} \frac{(r \ln(t))^n}{n!} \right) \mid rk>$$

$$= e^{(r \ln(t))} \mid rk> \; = \; t^r \mid rk> .$$

**Proposition 2.** $A^P \mid rk> = A^k \mid rk>$

**Proof:**

$$A^P \mid rk> = e^{\ln(A^P)} \mid rk>$$

$$= e^{(P \ln(A))} \mid rk> = \left( \sum_{n=0}^{\infty} \frac{(P \ln(A))^n}{n!} \right) \mid rk>$$

$$= \sum_{n=0}^{\infty} \frac{\ln(A)^n P^n \mid rk>}{n!} = \left( \sum_{n=0}^{\infty} \frac{\ln(A)^n k^n}{n!} \right) \mid rk>$$

$$= \left( \sum_{n=0}^{\infty} \frac{(k \ln(A))^n}{n!} \right) \mid rk> = e^{(k \ln(A))} \mid rk>$$

$$= A^k \mid rk>$$

**Proposition 3.**

$$U \mid rk> = t^r A^k \operatorname{rank}( \mathcal{H}^{k,r}(C) ) \mid rk>$$

**Proof:**

$$U \mid rk> = t^F A^P M \mid rk>$$

$$= t^F A^P [M \mid rk>] = t^F A^P [\operatorname{rank}( \mathcal{H}^{k,r}(C) ) \mid rk>]$$

$$= t^F [A^P \mid rk>] \operatorname{rank}( \mathcal{H}^{k,r}(C) ) = t^F [A^k \mid rk>] \operatorname{rank}( \mathcal{H}^{k,r}(C) )$$

$$= [t^F \mid rk>] A^k \operatorname{rank}( \mathcal{H}^{k,r}(C) ) = t^r A^k \operatorname{rank}( \mathcal{H}^{k,r}(C) ) \mid rk> .$$

With all these propositions, the quantum model of the Khovanov homology for tangles is given by

$$\mathrm{tr}(U) = \mathrm{tr}(t^F A^P M) = Tangle\ polynomial = \mathcal{P}(T)_{\varepsilon,Z} \tag{16}$$

then we have

$$\mathcal{P}(T)_{\varepsilon,Z} = \sum_r \left( \sum_k \langle rk \mid U \mid rk \rangle \right) \tag{17}$$

$$\mathcal{P}(T)_{\varepsilon,Z} = \sum_r \left( \sum_k t^r A^k \,\mathrm{rank}(\, \mathcal{H}^{k,r}(C)\, ) \langle rk \mid rk \rangle \right) \tag{18}$$

$$\mathcal{P}(T)_{\varepsilon,Z} = \sum_r \sum_k t^r A^k \,\mathrm{rank}(\, \mathcal{H}^{k,r}(C)\, ) \tag{19}$$

Finally we have that

$$\lim_{\tau \to 0} \mathrm{Tr}(t^F A^p M\, \mathrm{e}^{(-\tau H)}) = \lim_{\tau \to 0} \int\!\!\int \mathrm{e}^{(-S(\phi,\,\psi))}\, dBoson\ dFermion$$

$$= \mathcal{P}(T)_{\varepsilon,Z} = \sum_r \sum_k t^r A^k \,\mathrm{rank}(\, \mathcal{H}^{k,r}(C)\, ) \tag{20}$$

We consider (20) as a stringy topological quantum computer which is able to compute the tangle polynomial in the Khovanov homology for tangles introduced by Lauda. In the next section we will consider the particular case of a open-closed stringy topological quantum computer which is able to compute the relevant Frobenius algebras via cobordisms.

We conclude this section with a preliminary simulation of a simple open-closed string topological quantum computer based on qubits made up of non-Abelian string anyons with Majorana modes, resulting from adding topological lattice defects as genons to the Abelian phase of the Kitaev honeycomb model. We obtain a computation of the tangle polynomials for the following tangle :

$[\![T]\!] :=$ 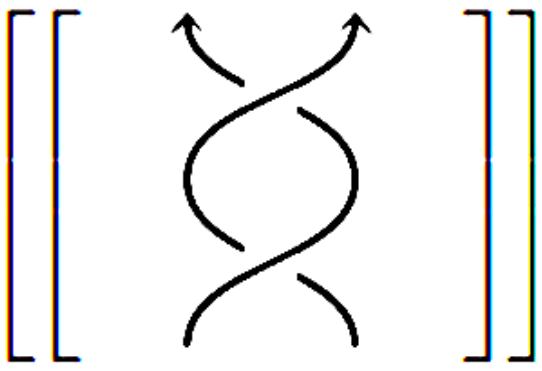 (21)

In the case of the open-closed stringy topological quantum computer of the Bar-Natan kind we have

$$\mathcal{H}^{2,0}(C') = \langle 1 \otimes y + y \otimes y \rangle, \ \mathcal{H}^{4,0}(C) = \langle 1 \otimes y + y \otimes 1 \rangle, \ \mathcal{H}^{8,2}(C) = \langle 1 \otimes x \otimes 1 \rangle, \ \mathcal{H}^{12,2}(C') = \langle 1 \otimes 1 \otimes 1 \rangle. \quad (22)$$

From (22) and using the proposition 3 we derive that

$$U \,|\, \Phi_{0,2} > = A^2 \,|\, \Phi_{0,2} > , \quad U \,|\, \Phi_{0,4} > = A^4 \,|\, \Phi_{0,4} > , \quad U \,|\, \Phi_{2,8} > = t^2 A^8 \,|\, \Phi_{0,4} > , \quad U \,|\, \Phi_{2,12} > = t^2 A^{12} \,|\, \Phi_{0,4} > \quad (23)$$

According with (17) and (22) we have that

$$\mathcal{P}_{Z_{\mathrm{BN}}}(T) = < \Phi_{0,2} \,|\, U \,|\, \Phi_{0,2} > + < \Phi_{0,4} \,|\, U \,|\, \Phi_{0,4} > + < \Phi_{2,8} \,|\, U \,|\, \Phi_{2,8} > + < \Phi_{2,12} \,|\, U \,|\, \Phi_{2,12} > \quad (24)$$

Using (23 ) and (24) we obtain

$$\mathcal{P}_{Z_{\mathrm{BN}}}(T) = A^2 + A^4 + t^2 A^8 + t^2 A^{12} \quad (25)$$

In the case of the open-closed stringy topological quantum computer of the Khovanov kind we have

$$\mathcal{H}^{2,0}(C') = \langle y \otimes y \rangle, \ \mathcal{H}^{4,0}(C) = \langle 1 \otimes y + y \otimes 1 \rangle, \ \mathcal{H}^{6,1}(C) = \left\langle \begin{pmatrix} 1 \otimes y \\ 1 \otimes y \end{pmatrix} \right\rangle,$$

$$\mathcal{H}^{8,1}(C) = \left\langle \begin{pmatrix} 1 \otimes 1 \\ 1 \otimes 1 \end{pmatrix} \right\rangle, \ \mathcal{H}^{8,2}(C) = \langle y \otimes 1 \otimes y \rangle, \ \mathcal{H}^{10,2}(C) = \langle 1 \otimes 1 \otimes y, y \otimes 1 \otimes 1 \rangle, \mathcal{H}^{12,2}(C) = \langle 1 \otimes 1 \otimes 1 \rangle, \quad (26)$$

From (26) and using the proposition 3 we derive that

$$U \,|\, \Phi_{0,2} > = A^2 \,|\, \Phi_{0,2} > , \quad U \,|\, \Phi_{0,4} > = A^4 \,|\, \Phi_{0,2} > , \quad U \,|\, \Phi_{1,6} > = t A^6 \,|\, \Phi_{1,6} > ,$$

$$U \,|\, \Phi_{1,8} > = t A^8 \,|\, \Phi_{1,8} > , \quad U \,|\, \Phi_{2,8} > = t^2 A^8 \,|\, \Phi_{2,8} > , \quad U \,|\, \Phi_{2,10} > = 2\, t^2 A^{10} \,|\, \Phi_{2,10} > ,$$

$$U \,|\, \Phi_{2,12} > = t^2 A^{12} \,|\, \Phi_{2,12} > . \quad (27)$$

According with (17) and (26) we have that

$$\mathcal{P}_{Z_{\mathrm{Kh}}}(T) := < \Phi_{0,2} \,|\, U \,|\, \Phi_{0,2} > + < \Phi_{0,4} \,|\, U \,|\, \Phi_{0,4} > + < \Phi_{1,6} \,|\, U \,|\, \Phi_{1,6} > + < \Phi_{1,8} \,|\, U \,|\, \Phi_{1,8} >$$
$$+ < \Phi_{2,8} \,|\, U \,|\, \Phi_{2,8} > + < \Phi_{2,10} \,|\, U \,|\, \Phi_{2,10} > + < \Phi_{2,12} \,|\, U \,|\, \Phi_{2,12} > \quad (28)$$

Using (27 ) and (28) we obtain

$$\mathcal{P}_{Z_{\mathrm{Kh}}}(T) = A^2 + A^4 + t\,A^6 + t\,A^8 + t^2\,A^8 + 2\,t^2\,A^{10} + t^2\,A^{12} \qquad (29)$$

## 4. USING LEO II IN STRING TOPOLOGICAL QUANTUM COMPUTING

The generators of the Loop braid group are [7]

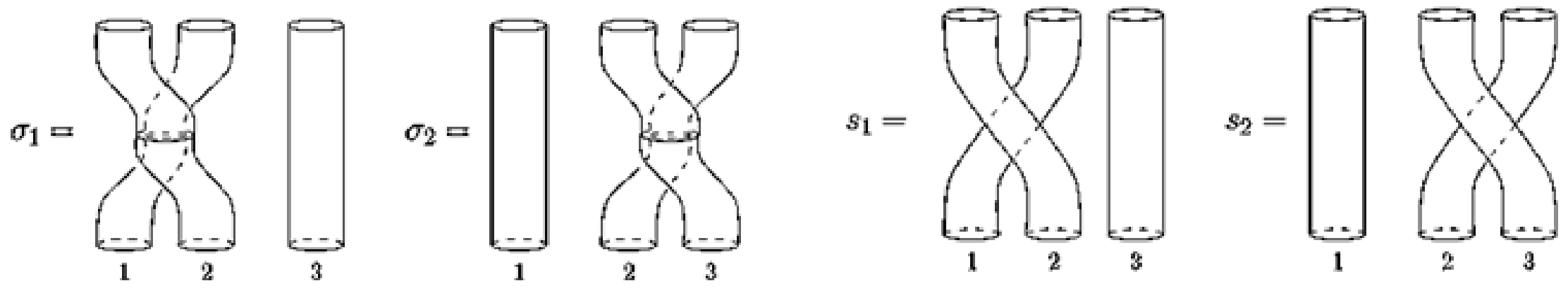

The generators of the open-closed cobordism category are [8,9]

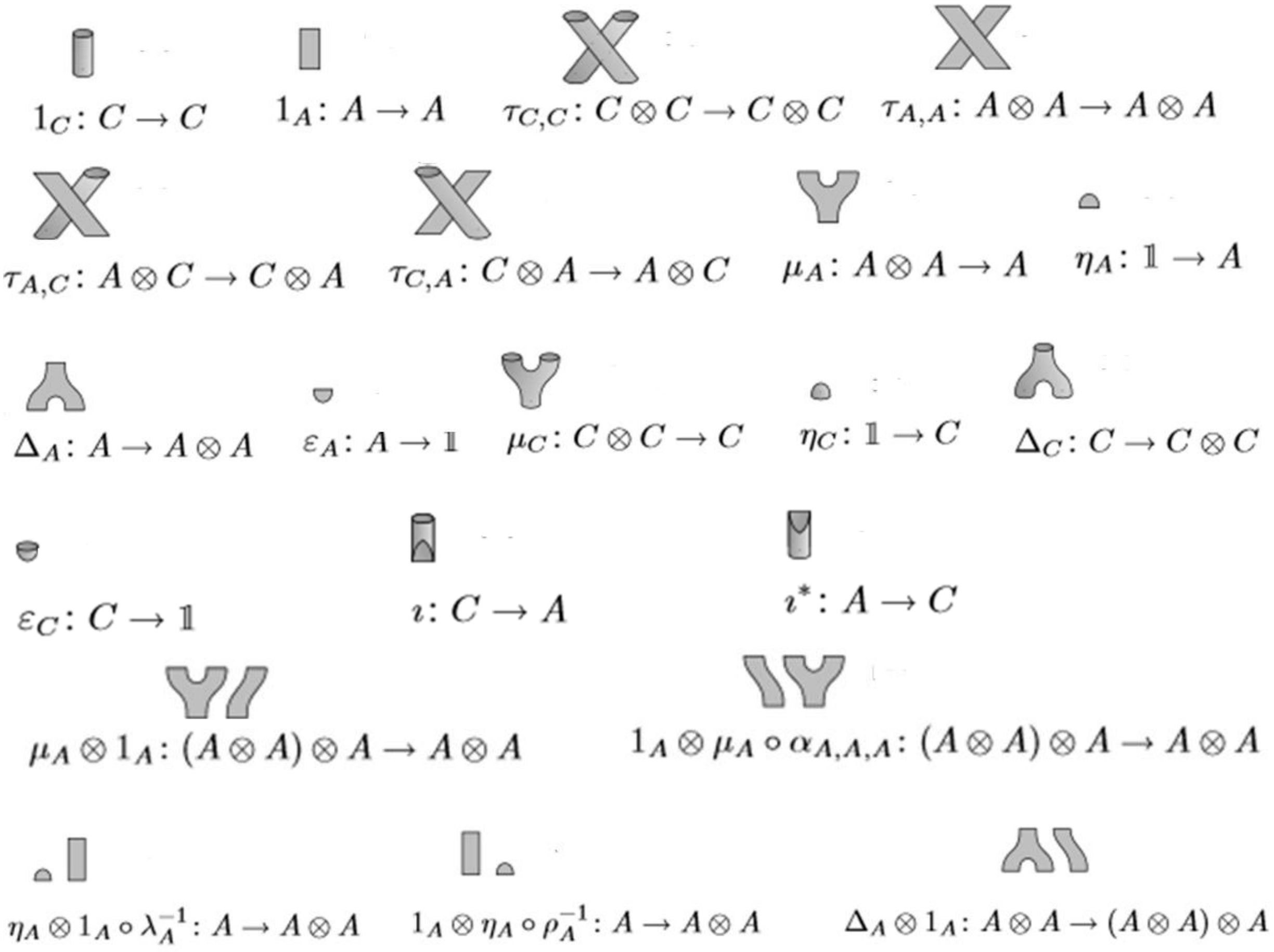

$\alpha^{-1}_{A,A,A} \circ 1_A \otimes \Delta_A \colon A \otimes A \to (A \otimes A) \otimes A$

$\lambda_A \circ \varepsilon_A \otimes 1_A \colon A \otimes A \to A$

$\rho_A \circ 1_A \otimes \varepsilon_A \colon A \otimes A \to A$

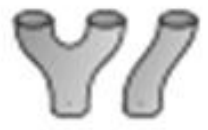

$\mu_C \otimes 1_C \colon (C \otimes C) \otimes C \to C \otimes C$

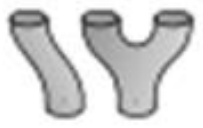

$1_C \otimes \mu_C \circ \alpha_{C,C,C} \colon (C \otimes C) \otimes C \to C \otimes C$

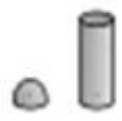 $\mapsto \quad \eta_C \otimes 1_C \circ \lambda_C^{-1} \colon C \to C \otimes C$

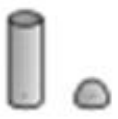 $\mapsto \quad 1_C \otimes \eta_C \circ \rho_C^{-1} \colon C \to C \otimes C$

$\Delta_C \otimes 1_C \colon C \otimes C \to (C \otimes C) \otimes C$

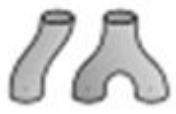 $\mapsto \quad \alpha^{-1}_{C,C,C} \circ 1_C \otimes \Delta_C \colon C \otimes C \to (C \otimes C) \otimes C$

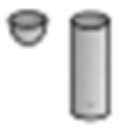 $\mapsto \quad \lambda_C \circ \varepsilon_C \otimes 1_C \colon C \otimes C \to C$

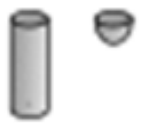 $\mapsto \quad \rho_C \circ 1_C \otimes \varepsilon_C \colon C \otimes C \to C$

where *C*is the Hilbert space for the quantum states of the closed strings and *A* is the Hilbert space for the quantum states of the open string. The generators of the open-closed cobordism category are considered as the basic Feynman diagrams for the interactions between the closed and open strings. At the same way, the generators of the Loop braid group are considered as the basic Feynman diagrams for the interactions between the closed strings. Combining the generators of the open-closed cobordism category and the generators of the Loop braid group it is possible to draw the more general Feynman diagrams for the interactions between open and closed strings. Two examples of such Feynman diagrams are showed at figure 1.

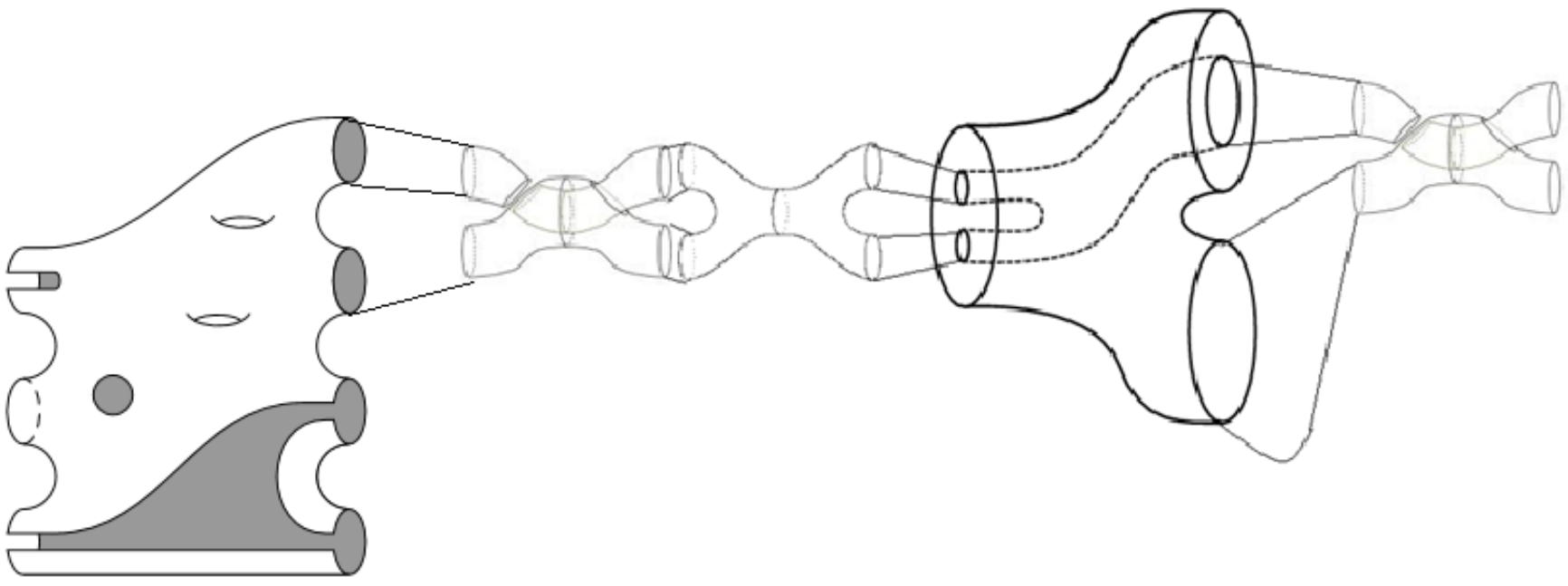

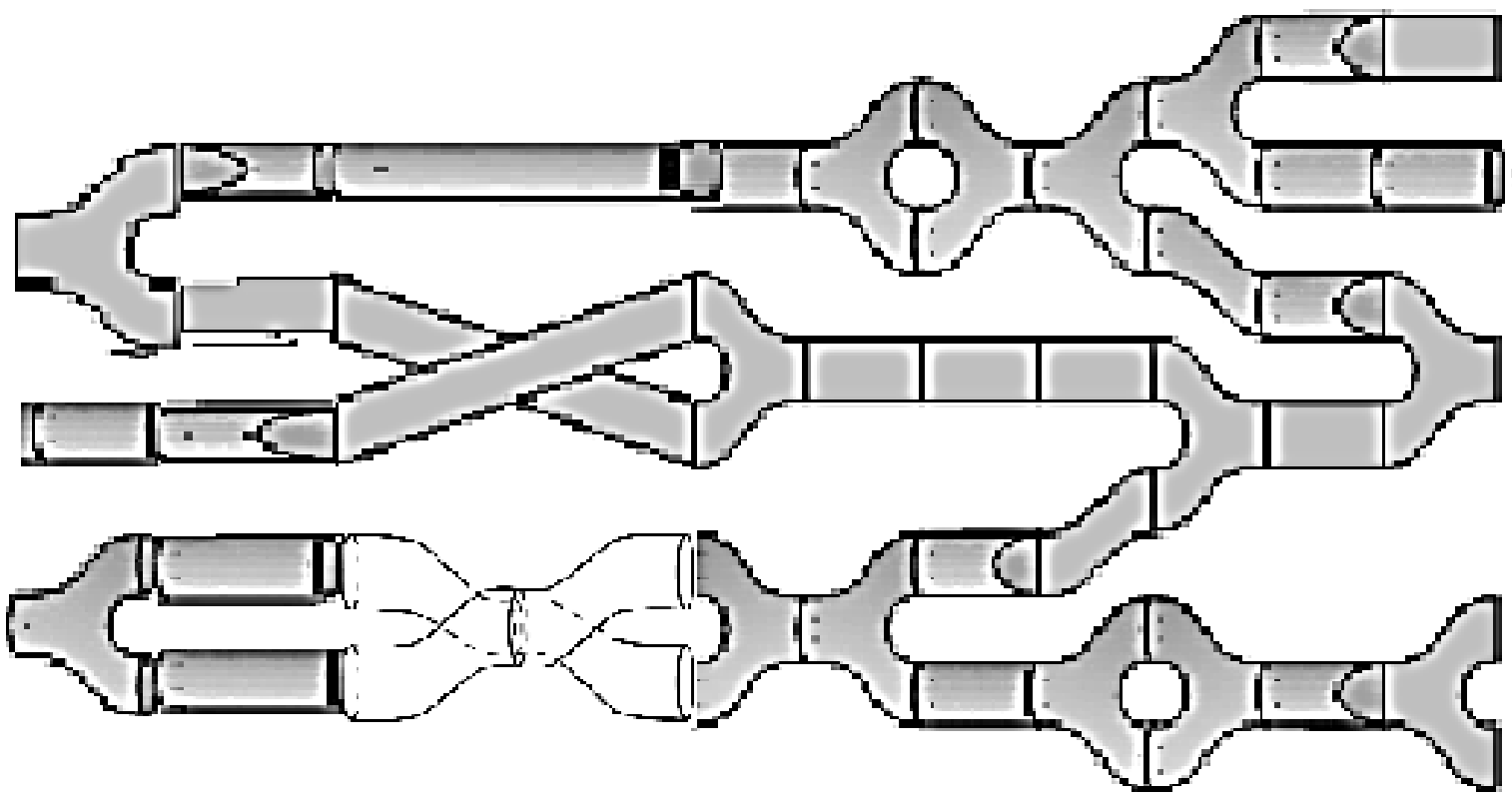

**Figure 1.** Two examples of Feynman diagrams for interactions between open and closed strings constructed from the generators of the open-closed cobordidms category and the generators of the Loop braid group.

The two Feynman diagrams showed at figure 1 can be considered respectively as representations of string topological quantum computations. A String Topological Quantum Computation is theoretically performed by braiding the paths of non-abelian string- anyons (sanyons), changing their states. The output of the computation uniquely depends on how the sanyons have been braided by the Loop Braid Group. Small disturbances do not unravel the loop braids, making the computation resistant to errors and decoherence. More in general, a String Topological Quantum Computation is theoretically performed by evolving non-abelian string- anyons (sanyons) using the Loop Braid Group and the open-closed cobordism category. A simple picture of a string topological quantum computer is showed at figure 2.

**Vacuum**

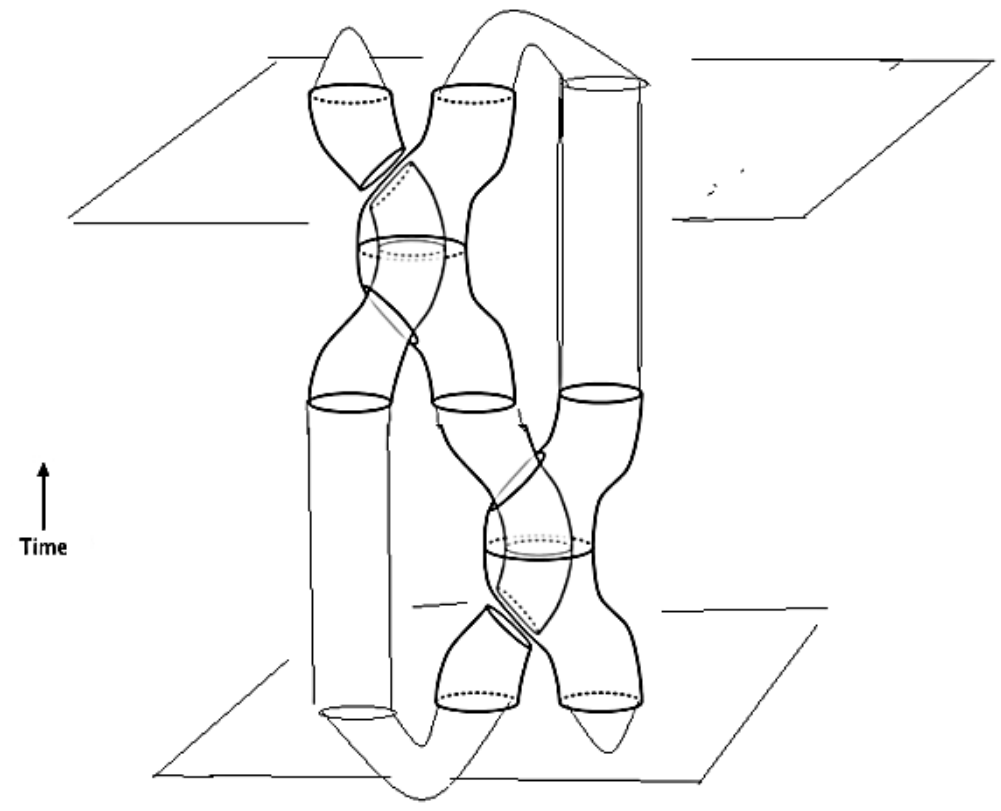

**Vacuum**

**Figure 2.** An illustration of a string topological quantum computer.

The mathematical design of an open-closed stringy topological quantum computer is involved with computations and theorem proving for generalized Frobenius algebras associated with open-closed cobordisms; and the Loop Braid Group. Such computations and theorem proving can be performed automatically using the ATP Leo II with the TPTP-thf language

An example of the application of ATP in the mathematical setup of an open-closed stringy topological quantum computer is as follows. Consider to prove the following theorem in Frobenius algebras:

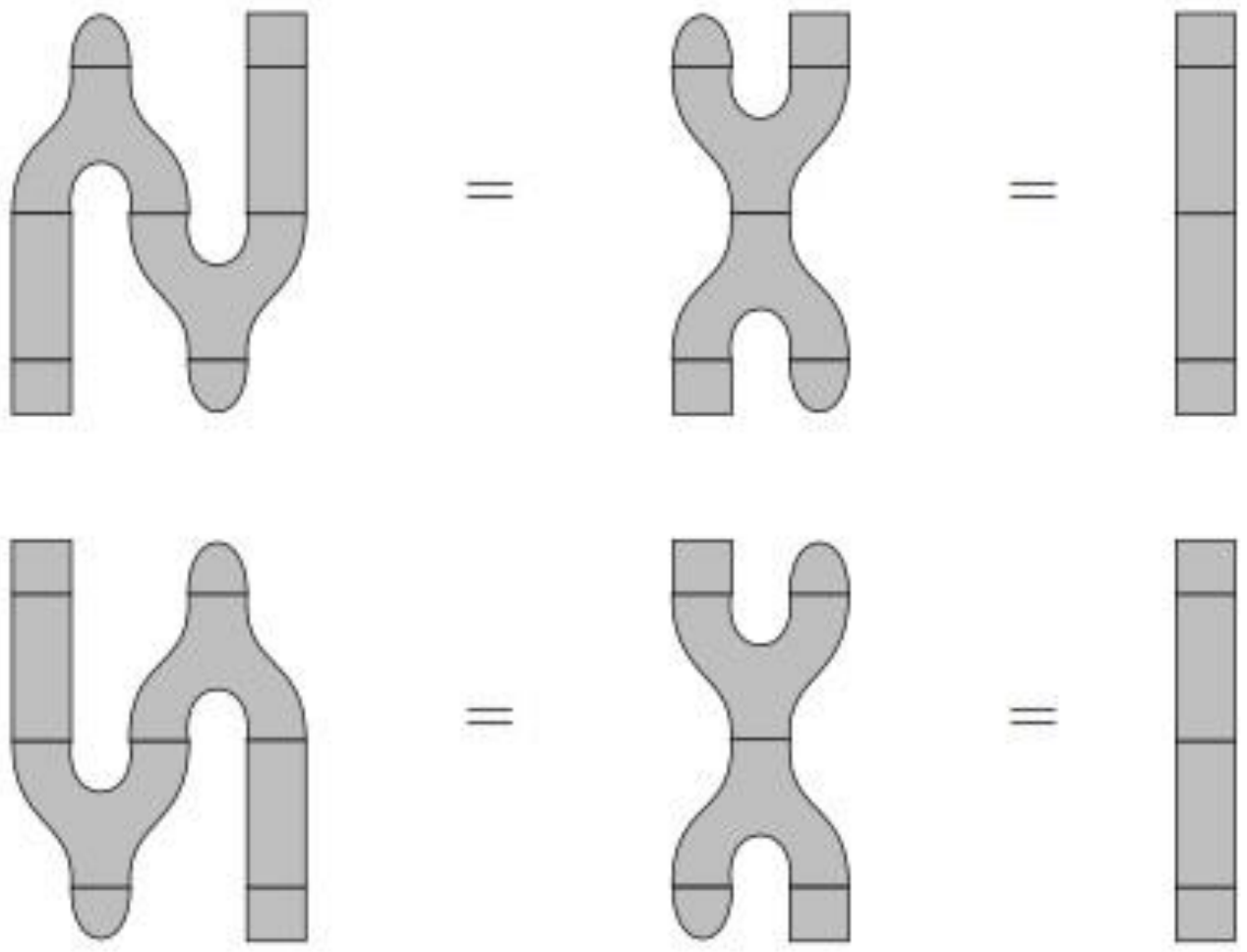

This theorem is proved using the ATP Leo II, using TPTP thf with the following code

```
thf(alpha_decl,type,(alpha: $aaxa > $axaa)).
thf(invalpha_decl,type,(invalpha: $axaa > $aaxa )).
thf(mu_decl,type,(mu: $aa > $a )).
thf(eta_decl,type,(eta: $i > $a )).
thf(muid_decl,type,(muid: $aaxa > $aa )).
thf(idmu_decl,type,(idmu: $axaa > $aa )).
thf(etaid_decl,type,(etaid: $ja > $aa )).
thf(ideta_decl,type,(ideta: $ai > $aa )).
thf(lamb_decl,type,(lamb: $ja > $a )).
thf(rho_decl,type,(rho: $ai > $a )).
thf(delta_decl,type,(delta: $a > $aa )).
thf(deltaid_decl,type,(deltaid: $aa > $aaxa )).
thf(iddelta_decl,type,(iddelta: $aa > $axaa )).
thf(epsilon_decl,type,(epsilon: $a > $i )).
thf(invlamb_decl,type,(invlamb: $a > $ja )).
thf(invrho_decl,type,(invrho: $a > $ai )).
thf(epsilonid_decl,type,(epsilonid: $aa > $ja )).
thf(idepsilon_decl,type,(idepsilon: $aa > $ai )).
thf(id_decl,type,(id: $a > $a )).
thf(beta1_decl,type,( beta1: $a > $ai )).
thf(beta2_decl,type,(beta2: $a > $ja )).
thf(invbeta1_decl,type,(invbeta1: $ai > $a )).
thf(invbeta2_decl,type,(invbeta2: $ja > $a )).
thf(axio1,axiom,(! [X: $ja] :
      ( (lamb @ X)  = (mu @ (etaid @ X ) )   )    )).
thf(axio2,axiom,(! [X: $ai] :
      ( (rho @ X)  = (mu @ (ideta @ X ) )   )    )).
thf(axio3,axiom,(! [X: $aaxa] :
      ( (mu @ (idmu @ (alpha @ X))   )  = (mu @ (muid @ X) )   )    )).
thf(axio4,axiom,(! [X: $a] :
      ( (iddelta @ (delta @ X) )  = (alpha @ (deltaid @ (delta @ X)) ) )   )).
thf(axio5,axiom,(! [X: $a] :
      (  (epsilonid @ (delta @ X) )  = (invlamb @ X )    )   )).
thf(axio6,axiom,(! [X: $a] :
```

```
      (  (idepsilon @ (delta @ X) )  = (invrho @ X )     )    )).
thf(axio7,axiom,(! [X: $aa] :
      (  (muid @ (invalpha @ (iddelta @ X))     )  = (delta @ (mu @ X))     )
)).
thf(axio8,axiom,(! [X: $aa] :
      (  (idmu @ (alpha @ (deltaid @ X))     )  = (delta @ (mu @ X))     )    )).
thf(axio9,axiom,(! [X: $a] :
      (  (mu @ (ideta @ (beta1 @ X))     )  = (id @ X)     )    )).
thf(axio10,axiom,(! [X: $a] :
      (  (mu @ (etaid @ (beta2 @ X))     )  = (id @ X)     )    )).
thf(axio11,axiom,(! [X: $a] :
      (  (invbeta1 @ (idepsilon @ (delta @ X))     )  = (id @ X)     )    )).
thf(axio12,axiom,(! [X: $a] :
      (  (invbeta2 @ (epsilonid @ (delta @ X))     )  = (id @ X)     )    )).
thf(axio13,axiom,(! [X: $a] :
      (  (id @ (id @ X)     )  = (id @ X)     )    )).
thf(axio14,axiom,(! [X: $a] :
      (  (delta @ (id @ X)     )  = (delta @ X)     )    )).
thf(conj,conjecture,(! [X: $ja] :
      ( (idepsilon @ (idmu @ (alpha @ (deltaid @ (etaid @ X))) )  ) =
        (idepsilon @ (delta @ (mu @ (etaid @ X))))
                )  )).
thf(conj2,conjecture,(! [X: $a] :
      ( (invbeta1 @ (idepsilon @ (delta @ (mu @ (etaid @ (beta2 @ X)))) )  ) =
        (id @ X)
                )  )).
thf(conj3,conjecture,(! [X: $ai] :
      ( (epsilonid @ (muid @ (invalpha @ (iddelta @ (ideta @ X))) )  ) =
        (epsilonid @ (delta @ (mu @ (ideta @ X))))
                )  )).
thf(conj4,conjecture,(! [X: $a] :
      ( (invbeta2 @ (epsilonid @ (delta @ (mu @ (ideta @ (beta1 @ X)))) )  ) =
        (id @ X)
                )  )).
```

and the corresponding output from Leo II is:

```
%**** Beginning of derivation protocol ****
% SZS output start CNFRefutation

% SZS output end CNFRefutation

%**** End of derivation protocol ****
%**** no. of clauses in derivation: 163 ****
%**** clause counter: 162 ****

% SZS status Theorem for /tmp/SystemOnTPTPFormReply39420/SOT_ZN9MIY :
(rf:0,axioms:13,ps:3,u:6,ude:true,rLeibEQ:true,rAndEQ:true,use_choice:true,use_e
xtuni:true,use_extcnf_combined:true,expand_extuni:false,foatp:e,atp_timeout:7,at
p_calls_frequency:10,ordering:none,proof_output:1,protocol_output:false,clause_c
ount:162,loop_count:0,foatp_calls:1,translation:fof_full)

% END OF SYSTEM OUTPUT
```

```
RESULT: SOT_ZN9MIY - LEO-II---1.6.2 says Theorem - CPU = 0.12 WC = 0.17
OUTPUT: SOT_ZN9MIY - LEO-II---1.6.2 says CNFRefutation - CPU = 0.12 WC = 0.18
```

Other example. To prove the following theorem in Frobenius algebras (zig-zag identities).

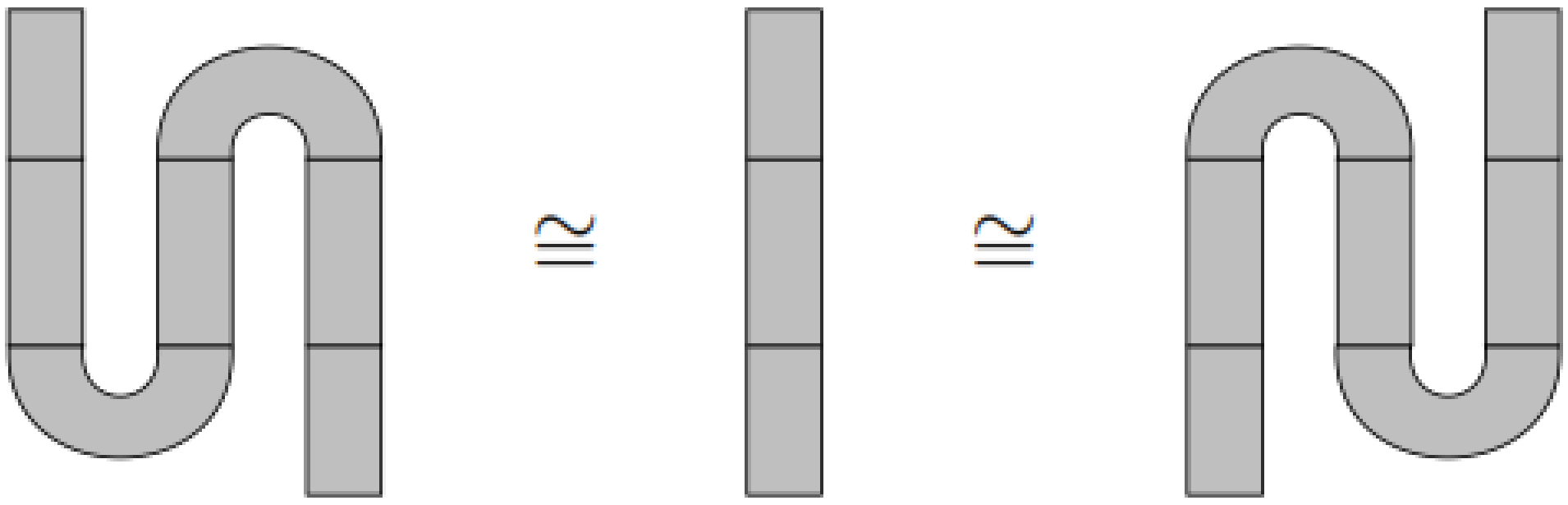

where the open pairing and open copairing are defined respectively by the following open-closed cobordisms

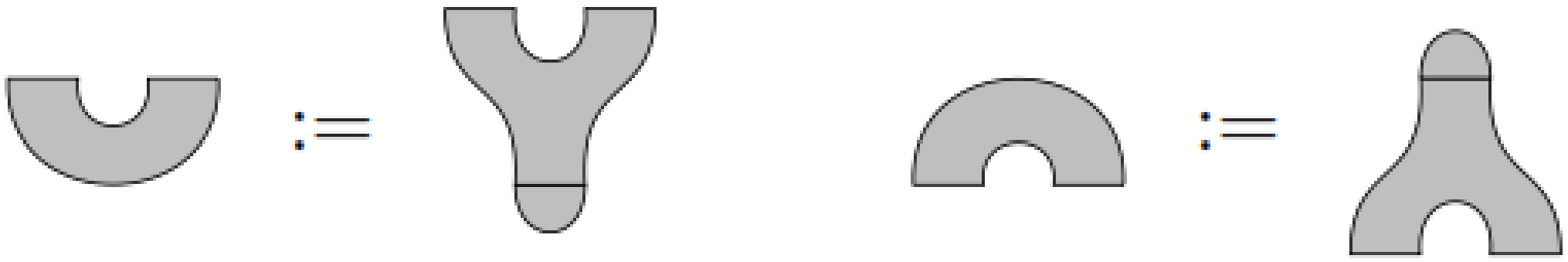

The proof of the open zig-zag identities is performed automatically using the ATP Leo II with TPTP thf language with the following code

```
thf(conj5,conjecture,(! [X: $a] :
      ( (invbeta2 @ (epsilonid @ (muid @ (invalpha @ (iddelta @ (ideta @ (beta1
@ X)))) ))  ) =
        (id @ X)
                )  )).

thf(conj6,conjecture,(! [X: $a] :
      ( (invbeta1 @ (idepsilon @ (idmu @ (alpha @ (deltaid @ (etaid @ (beta2 @
X)))) ))  ) =
        (id @ X)
                )  )).

thf(conj7,conjecture,(! [X: $a] :
      ((invbeta2 @ (epsilonid @ (muid @ (invalpha @ (iddelta @ (ideta @ (beta1 @
X)))) ))  ) =
        (invbeta1 @ (idepsilon @ (idmu @ (alpha @ (deltaid @ (etaid @ (beta2 @
X)))) ))  )                   )  )).
```

and the corresponding output from Leo II is

```
% END OF SYSTEM OUTPUT
% RESULT: SOT_DDEXiT - LEO-II---1.6.2 says Theorem - CPU = 0.08 WC = 0.13
```

```
% OUTPUT: SOT_DDEXiT - LEO-II---1.6.2 says CNFRefutation - CPU = 0.08 WC = 0.14
```

Other example. To prove the following theorem in Frobenius algebra

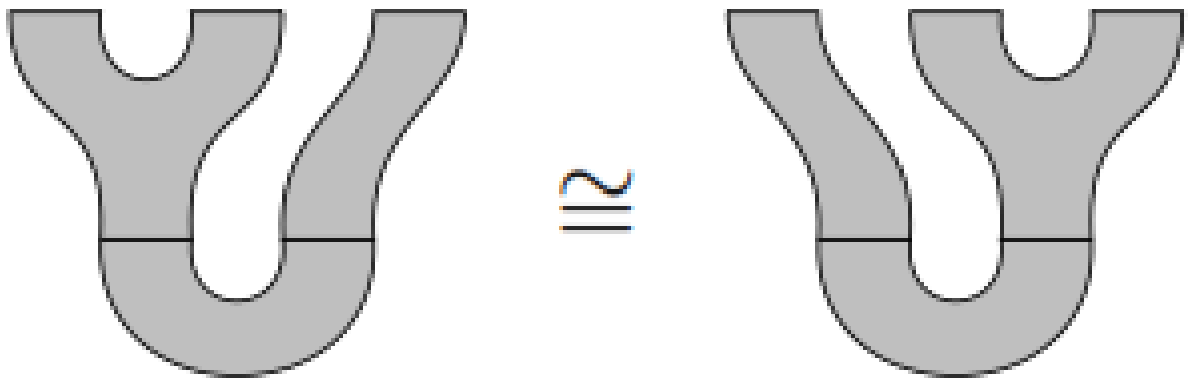

We use the ATP Leo II with TPTP thf language with the following code

```
thf(conj8,conjecture,(! [X: $aaxa] :
      ( (epsilon @ (mu @ (muid @ X))  ) =
        (epsilon @ (mu @ (idmu @ (alpha @ X))))
                )  )).
```

and the corresponding output is

```
% END OF SYSTEM OUTPUT
% RESULT: SOT_o1E2T7 - LEO-II---1.6.2 says Theorem - CPU = 0.00 WC = 0.04
% OUTPUT: SOT_o1E2T7 - LEO-II---1.6.2 says CNFRefutation - CPU = 0.00 WC = 0.04
```

The corresponding theorem for the copairing is

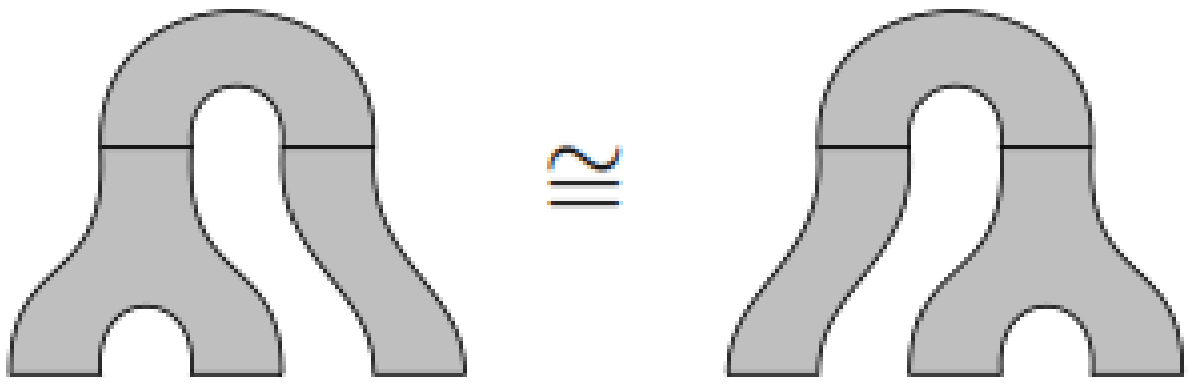

and using the ATP Leo II with TPTP thf language with the following code

```
thf(conj9,conjecture,(! [X: $i] :
      ( (alpha @ (deltaid @ (delta @ (eta @ X))  )) =
         (iddelta @ (delta @ (eta @ X)))
                )  )).
```

The output from Leo II is

```
% END OF SYSTEM OUTPUT
% RESULT: SOT_K7V_sv - LEO-II---1.6.2 says Theorem - CPU = 0.02 WC = 0.05
% OUTPUT: SOT_K7V_sv - LEO-II---1.6.2 says CNFRefutation - CPU = 0.02 WC = 0.05
```

Other example. Prove the following theorem in Frobenius algebra

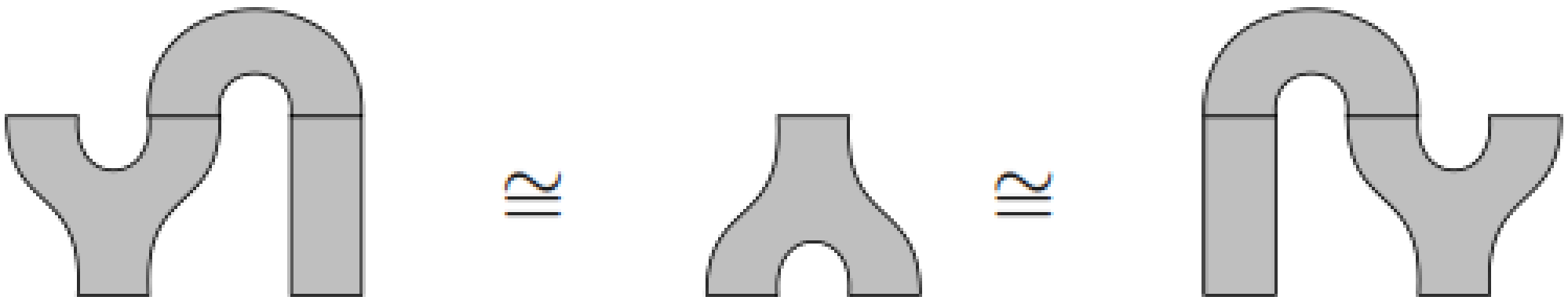

This theorem is proved automatically using the ATP Leo II with the following TPTP thf code

```
thf(conj10,conjecture,(! [X: $a] :
      ( (muid @ (invalpha @ (iddelta @ (ideta @ (beta1 @ X))  ))  ) =
        (delta @ X)
                 )  )).


thf(conj11,conjecture,(! [X: $a] :
      ( (idmu @ (alpha @ (deltaid @ (etaid @ (beta2 @ X))  ))  ) =
        (delta @ X)
                 )  )).

thf(conj12,conjecture,(! [X: $a] :
      ( (muid @ (invalpha @ (iddelta @ (ideta @ (beta1 @ X))  ))  ) =
        (idmu @ (alpha @ (deltaid @ (etaid @ (beta2 @ X))  ))  )
                 )  )).
```

and the corresponding output from Leo II is

```
% END OF SYSTEM OUTPUT
% RESULT: SOT_BR9hhI - LEO-II---1.6.2 says Theorem - CPU = 0.08 WC = 0.12
% OUTPUT: SOT_BR9hhI - LEO-II---1.6.2 says CNFRefutation - CPU = 0.08 WC = 0.13
```

Other example. The following theorem in Frobenius algebra

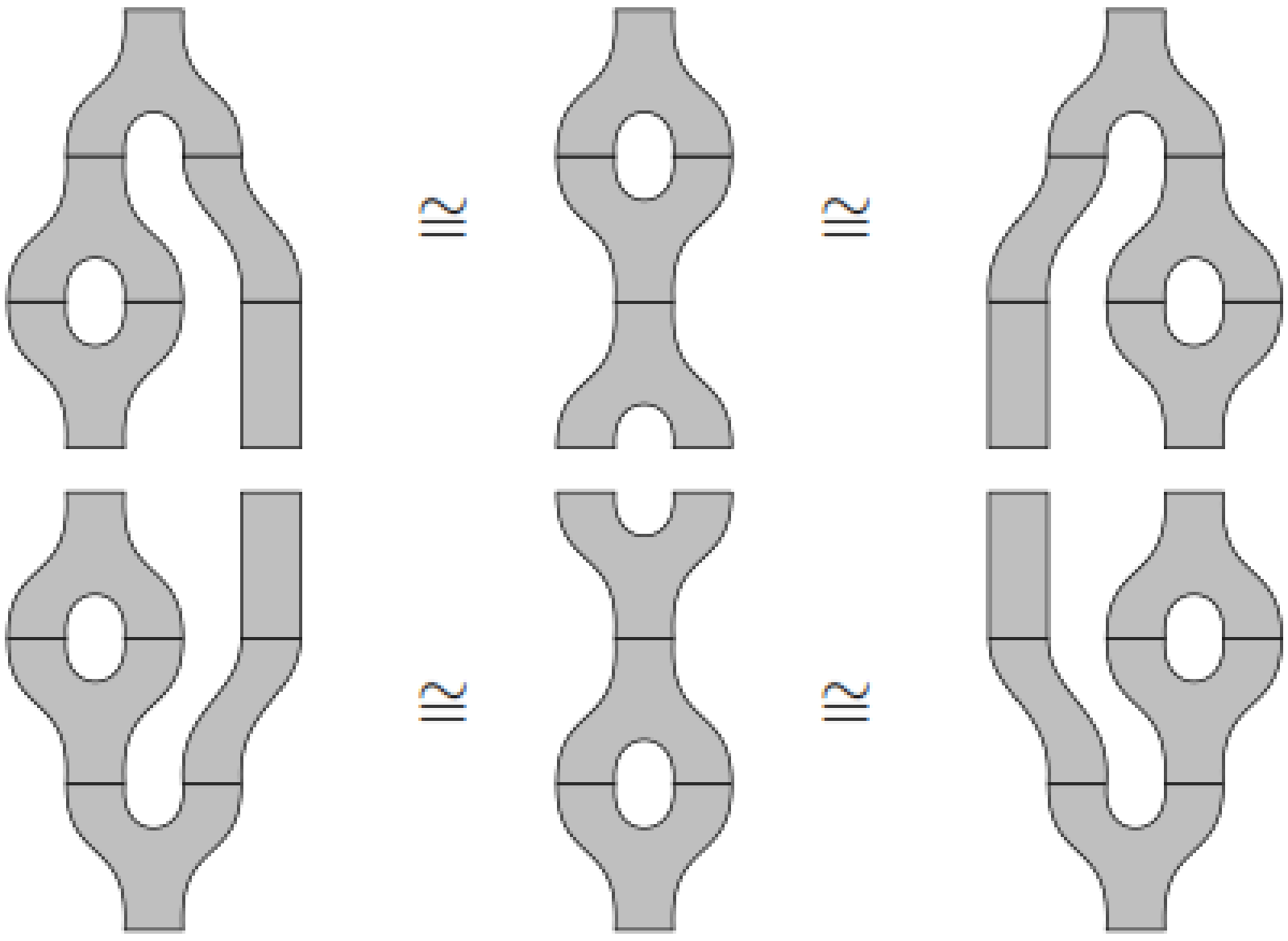

The first identity is proved automatically using the ATP Leo II with the following TPTP thf code

```
thf(conj13,conjecture,(! [X: $a] :
      ( (muid @ (invalpha @ (alpha @ (deltaid @ (delta @ X))))) =
        ( (delta @ (mu @ (delta @ X)) ))
                )  )).

thf(conj14,conjecture,(! [X: $a] :
      ( (muid @ (invalpha @ (alpha @ (deltaid @ (delta @ X))))) =
        (idmu @ (iddelta @ ((delta @ X)) ))
                )  )).

thf(conj15,conjecture,(! [X: $a] :
      ( (idmu @ (iddelta @ ((delta @ X)) )) =
        ( (delta @ (mu @ (delta @ X)) ))
                )  )).
```

and the corresponding output from Leo II is

```
END OF SYSTEM OUTPUT
% RESULT: SOT_qXxzea - LEO-II---1.6.2 says Theorem - CPU = 0.08 WC = 0.14
% OUTPUT: SOT_qXxzea - LEO-II---1.6.2 says CNFRefutation - CPU = 0.08 WC = 0.15
```

The second identity is proved automatically using the ATP Leo II with the following TPTP thf code

```
thf(conj16,conjecture,(! [X: $aa] :
      ( (mu @ (idmu @ (alpha @ (invalpha @ (iddelta @ X))))) =
        (mu @ (muid @ ((deltaid @ X)) ))
                )  )).

thf(conj17,conjecture,(! [X: $aa] :
      ( (mu @ (muid @ ((deltaid @ X)) )) =
        ( (mu @ (delta @ (mu @ X)) )) )  )).
```

```
thf(conj18,conjecture,(! [X: $aa] :
      ( (mu @ (idmu @ (alpha @ (invalpha @ (iddelta @ X))))) =
        (mu @ (delta @ ((mu @ X)) ))
                )  )).
```

and the corresponding output from Leo II is

```
% END OF SYSTEM OUTPUT
% RESULT: SOT_Ft33gx - LEO-II---1.6.2 says Theorem - CPU = 0.08 WC = 0.12
% OUTPUT: SOT_Ft33gx - LEO-II---1.6.2 says CNFRefutation - CPU = 0.08 WC = 0.13
```

## 5. CONCLUSIONS

A quantum model of Knovanov Homology for tangles, was presented using the supersymmetric string theory with open and closed strings. A generalized character valued index in superstring theory was proposed as an open-closed stringy topological quantum computer which is able to compute the tangle polynomial for an arbitrary tangle. Such open-closed stringy topological quantum computer was mathematically designed using the automatic theorem prover named Leo II which was able to prove theorems for the generalized Frobenius algebras associated with the Khovanov homology for tangles. Some simulations concerning with the computation of the tangle polynomials for simple tangles were performed using computer algebra software. The notion of an open-closed cobordism considered as a quantum gate was formulated using the generators of the Loop Braid group and the generators of the open-closed cobordism category. As a line for future research it is very interesting to investigate the application of the system named *Sledgehammer* to the verification of the proofs generated by Leo II in the case of the generalized Frobenius algebras involved in the mathematical design of an open-closed stringy topological quantum computer; and the possible applications of such topological quantum computer as a quantum automated theorem prover (QATP) based on quantum Hoare logic.